\documentclass[12pt,twosided]{article}
\usepackage{fullpage}
\usepackage{amsmath}
\usepackage{amsfonts}

\usepackage{amsthm}
\usepackage{bbm}

\usepackage{graphicx}
\usepackage{subfig}
\usepackage{epstopdf}
\usepackage{cite}
\usepackage{authblk}

\newtheorem{example}{Example}
\newtheorem{casestudy}{Case Study}

\newcommand{\etal}{\textit{et al.}}

\title{Designing Networks: A Mixed-Integer Linear Optimization Approach}

\author[1]{Chrysanthos E. Gounaris\thanks{Corresponding Author: gounaris@cmu.edu}}
\author[2]{Karthikeyan Rajendran}
\author[2,3]{Ioannis G. Kevrekidis}
\author[4,5]{Christodoulos A. Floudas}
\affil[1]{\small Department of Chemical Engineering, Carnegie Mellon University, Pittsburgh, PA 15213, USA}
\affil[2]{Department of Chemical and Biological Engineering, Princeton University, Princeton, NJ 08544, USA}
\affil[3]{Program in Applied and Computational Mathematics, Princeton University, Princeton, NJ 08544, USA}
\affil[4]{Texas A\&M Energy Institute, College Station, TX 77843, USA}
\affil[5]{Department of Chemical Engineering, Texas A\&M University, College Station, TX 77843, USA}

\date{}

\begin{document}

\maketitle

\newpage
\begin{abstract}
\normalsize
Designing networks with specified collective properties is useful in a variety of application areas, 
enabling the study of how given properties affect the behavior of network models, 
the downscaling of empirical networks to workable sizes, 
and the analysis of network evolution.
Despite the importance of the task, there currently exists a gap in our ability to systematically generate networks 
that adhere to theoretical guarantees for the given property specifications.
In this paper, we propose the use of Mixed-Integer Linear Optimization modeling and solution 
methodologies to address this {\em Network Generation Problem}.
We present a number of useful modeling techniques and apply them to mathematically express and 
constrain network properties in the context of an optimization formulation.
We then develop complete formulations for the generation of networks that attain specified levels of connectivity, spread, assortativity and robustness, and we illustrate these via a number of computational case studies.
\\

\end{abstract}

\newpage

\section{Introduction}\label{intro}
Detailed modeling and simulation of large, complex networks are increasingly being used in several fields (from biology, epidemiology and economics, to traffic management and supply-chains).  
As a result, there have been several empirical studies of the interaction between the structure of such networks and their dynamics (a ``network'' analog of
structure-property relations in materials science).
Detecting, understanding, and eventually exploiting this structure-dynamics relation in networks is thus a subject of intense interdisciplinary research.
Clearly, the construction of networks with prescribed structure is an integral part of this effort, and the purpose of this paper is
to present and illustrate a computer-assisted, mixed integer linear programming (MILP) approach to this important task.

There are many important examples where the topology of real world networks has been analyzed in great detail in order to
understand their formation and dynamics - the Internet holds probably pride of place among these.
A brief history of the development of Internet ``topology generators'' (the construction of networks that hopefully
model the Internet) can be found in Tangmunarunkit~\etal~\cite{Tang02network}.
The earliest internet topology generators were based on the \emph{Erd\"{o}s-R\'{e}nyi} model.  
Realistic simulations of dynamical models of the internet, however, required topology generators that more accurately represented its structure.  Hence, better models that take into account the hierarchy in the internet were later proposed~\cite{Calv97modeling}, and these were widely accepted until the breakthrough results by Faloutsos~\etal~\cite{Falo99power-law} revealed that the topology of the internet exhibits power-law dependence in terms of the degrees of the nodes.  
Networks with such power-law degree distributions are called \emph{scale-free} networks, and are in stark contrast with networks featuring the Poisson distribution of degrees created by topology generators based on the Erd\"{o}s-R\'{e}nyi model.  
This suggests that local properties of real networks are not accurately reproduced by simple graph models.
While a lot of attention has focused on designing realistic models of the Internet, the utility of graph models that produce power-law degree distributions is not restricted to the Internet domain alone.  
Several systems have been found to have power law degree distributions with different degree exponents. 
Newman~\cite{Newm03structure} and Albert and Barab\'{a}si~\cite{Albe02statistical} review graphs from numerous systems that adhere to a power-law degree distribution.  
These discoveries create the need for better alternative approaches to generating models of real world networks.
In response to this need, Barab\'{a}si and Albert~\cite{Bara99emergence} themselves proposed an evolutionary process to create scale-free networks. 
In their model, the network grows by the addition of nodes which preferentially attach to the existing nodes that are already well connected. This model is termed in the literature as the ``preferential attachment'' model.  
A different approach is used by Aiello~\etal~\cite{Aiel00random} to create random graphs with power-law degree distribution.
They first generate sequences of degrees from a power law distribution and then create networks that satisfy this degree distribution.

These two approaches exemplify two distinct (yet at some level related) approaches by which realistic network properties (i.e., the scale-free structure in this case) 
can be ``built in'' model networks.   
The first approach, which one might name \emph{complex network growth models} is to use {\em dynamical construction procedures} (evolutionary or growth processes) that give rise to networks with required properties.  
Examples of such models that have been extensively used in complex systems research are the Watts-Strogatz model for small-world networks~\cite{Watt98collective} and the aforementioned preferential attachment model that creates scale-free networks~\cite{Bara99emergence}.

While such complex network growth models may provide insights into the microscopic processes that govern the formation of real world networks,
they will in general create graphs that do not fully capture several properties of interest.
For example, Bu~\etal~\cite{Bu02distinguishing} show that Internet topology generators produce graphs which have different values of characteristic path lengths and clustering coefficients compared to the actual internet graph.
The second approach to network generation is built around the requirement to create graphs with specified properties. 
Such methods might be named \emph{network generation algorithms}.
The main feature of this second approach is the fact that {\em combinations of several properties} (beyond, for example, just the degree distribution) 
can be explicitly taken into account.

Newman~\cite{Newm03structure} provides a list of different graph properties that can be ``built in'' network generation algorithms,
and discusses their relevance in real systems. 
For instance, the average path length of social networks is found to be small for even very large networks, which characterizes what is known as the \emph{small-world property}. This has direct implications on the diffusion of information or the spread of epidemics on such networks. 
Shirley and Rushton~\cite{Shir05impacts} study the impact of the structure of four different types of network topologies, namely Erd\"{o}s-R\'{e}nyi, regular lattices, small-world, and scale-free, on epidemic dynamics. 
They show that in scale-free networks with very small clustering and short path lengths, the spread of the disease is faster, while in regular lattices with high clustering coefficients and long path lengths, the spread is the slowest among all other types of networks. 
Badham and Stocker~\cite{Badh10impact} report that increases in assortativity and
clustering coefficient are associated with smaller epidemics.
The effect of correlations between neighboring node degrees in a network has also been found to be a crucial property in determining system behavior in a variety of physical systems. 
Miao {\em et al.}~\cite{Miao08effects} study the effect of degree-degree correlations on the pinning controllability of networks and find that dissortative mixing enhances controllability. The results of Pechenick {\em et al.}~\cite{Pech12influence} show that chaotic
gene regulatory networks exhibit increasingly robust phenotypes with increasing assortativity. 
Pastor-Satorras {\em et al.}~\cite{Past01dynamical} report that disassortative networks enhance transport through them.  
They find that ``the pressure congestion for disassortative (assortative) networks is lower (bigger) than the one for uncorrelated networks.''
Watts and Strogatz~\cite{Watt98collective} report that oscillators connected by a small-world structure synchronize more readily than those connected by other network topologies. 
Hong~\etal~\cite{Hong04factors} find other important properties that affect the synchronizability of oscillators connected by a Watts-Strogatz network~\cite{Watt98collective}. 
They find that ``better synchronizability for the WS small-world network is induced as the heterogeneity
of the degree distribution is increased, as the characteristic path length is decreased, as the maximum betweenness centrality is reduced, or as the maximum degree is raised.''  
Thus, we find that properties like the maximum degree and the maximum betweenness centrality are also relevant for some dynamical network systems.
Ravasz and Barab\'{a}si~\cite{Rava03hierarchical} find that the clustering coefficients of the nodes are proportional to the inverse of the nodal degrees for four real networks (actor network, semantic web, world wide web and the internet). 
They propose this to be related to the hierarchical organization of the network. 
They suggest that ``the fundamental discrepancy between models and empirical measurements is rooted in a previously disregarded, yet generic feature of many real networks: their hierarchical topology.'' 
Thus, it may be important to maintain the ``degree-dependent clustering''  in models of network topology.

All these results clearly point to the necessity of being able to generate graphs with specified properties, not only to replicate realistic structures of real networks but also to create artificial networks that will result in desired functions. 
One of the main difficulties in accurately studying the effect of network properties on the function of networked systems is disassociating the effect of one property from the effect of others. 
Wrong conclusions will be derived when the network topology generator induces artificial correlations between different sets of properties.  
It is clear that graph-generating algorithms must be able to create networks with controllably-tunable combinations of graph properties.

A majority of existing network generation algorithms focus on constructing networks with specified degree statistics.
The Havel-Hakimi algorithm~\cite{Haki62realizability,Have55remark} is a deterministic algorithm that guarantees a network for any given graphical sequence of degrees. However, the algorithm always creates the same network for a given input degree sequence, which may be a limitation in certain applications. 
The ``pairing'' or ``configuration'' model~\cite{Boll80probabilistic,Worm80some} is a stochastic algorithm that alleviates this problem, but it suffers from non-uniform sampling.  
Strategies to generate uniformly distributed random graphs with specified degree sequences are discussed in a number of publications~\cite{Blit06sequential,Mil03uniform,Vige05efficient}. 
 The Chung-Lu algorithm~\cite{Chun02connected,Mill11efficient} is another probabilistic algorithm that creates a random graph from an ensemble of graphs with a given {\em expected} degree sequence.

Besides the degree sequence, clustering and assortativity are two other properties that have received considerable attention in the context of network topology generation algorithms.
Algorithms that tune clustering coefficients in scale-free networks are given in~\cite{Fu02evolving,Holm02growing}, while an algorithm for small-world networks can be found in~\cite{Wang08evolving}.
Several algorithms have also been proposed to construct graphs with specified clustering coefficients in conjunction with given degree distributions~\cite{Bans09exploring,Deij09random,Kim04performance,Volz04random}.
Procedures to create networks with specified assortativity have also been introduced~\cite{Guo06growing, Newm02assortative, Newm03mixing}.  
An algorithm that controls both the degree distribution and assortativity can be found in Xulvi-Brunet and Sokolov~\cite{Xulv04reshuffling}.  The problem of constructing graphs with given maximum degree and diameter that has the largest number of vertices has also been studied~\cite{Aleg86some,Berm84large,Fell95large}.
In addition, Nakano~\etal~\cite{nakano} discuss an algorithm to construct weighted graphs such that the total weight of the edges incident to a node is at least the given weight of the node.
Finally, other property combinations that have been addressed in the context of network generation algorithms include the
(a) degree distribution and degree-dependent clustering~\cite{OPTLpaper,BTER_2014,Serr05tuning,BTER_2012},
(b) joint degree-degree distribution~\cite{Doro02how,Webe07generation},
(c) joint degree-degree distribution and modularity~\cite{Zhou10generation},
(d) degree-degree correlations and degree-dependent clustering~\cite{Pusc08generating}, and
(e) eigenvector centrality~\cite{Nico12controlling}.

The majority of the above algorithms cannot guarantee the generation of networks that satisfy exactly the properties specifications, or that lie within bounded deviations from the desired values.
In our earlier work~\cite{OPTLpaper}, we introduced a method to specify degree distribution in conjunction with a specified number of triangles among nodes of given degree (i.e., degree-dependent clustering levels).
The method provided theoretical guarantees that a network satisfying the imposed specification would be constructed, if such a network existed.
It also constituted the first time in the open literature that a network generation approach based on mathematical optimization was proposed.  In particular, the problem was posed as a Mixed-Integer Linear Problem (MILP), and as a result, the method offered a number of advantages over traditional approaches.  
For example, in addition to identifying a network that possesses a collection of desired properties, an MILP-based framework is typically suitable to
(a) certify if no such network exists,
(b) identify multiple, or even all, non-isomorphic such networks,
(c) identify the network that is most similar --under a suitable metric-- to a reference network,
(d) identify the network that minimizes/maximizes some application-specific objective,
(e) determine the smallest size $N_{0}$ for a network to admit all properties, and
(f) identify networks that approximately satisfy the properties within given tolerances.

Posing the problem as an MILP allows for employing standard, well-developed solution algorithms (e.g., branch-and-bound, branch-and-cut) to deal with the complex decision-making that is involved.
It is thus a promising approach towards developing the capability to systematically design networks that possess desired properties with explicit guarantees of achieving the specifications exactly or within bounded deviations.
In this paper, we demonstrate the scope and generality of the MILP-based approach, and in particular the associated modeling techniques required to handle a wide variety of properties.
The remainder of the paper is structured as follows: after formalizing the Network Generation Problem, Section~\ref{modeling} introduces the modeling techniques that can be used to formulate the mathematical expression and specification of any network property of interest in an optimization context.
Section~\ref{popular} then develops model ``snippets'' that correspond to a number of popular network properties.
Section~\ref{enhancements} addresses modeling and implementation issues that can enhance the performance of optimization algorithms in the context of the Network Generation Problem.
Section~\ref{studies} discusses how to put together complete optimization formulations for the design of networks that have tunable levels of connectivity, spread, assortativity and robustness.  
Although by no means exhaustive for the capabilities of the proposed optimization-based framework, these illustrative applications were selected so as to showcase its versatility and flexibility with regards to network properties that can be addressed.  
A number of computational case studies are also provided for illustration.  
Finally, Section~\ref{conclusions} presents our concluding remarks.
Throughout the paper, we use boldface to indicate vector quantities.

\section{Modeling the Network Generation Problem}\label{modeling}
The Network Generation Problem (NGP) calls for the identification of the connections (edges) that need to exist between its nodes in order for the network to attain specific values, or values within specific ranges, for a set of properties of interest.
By ``properties'' we refer to various characteristic quantities of networks including global metrics (e.g., the global clustering coefficient), local metrics (e.g., a node's local clustering coefficient), nodal distributions (e.g., the degree sequence), network distances (e.g., shortest paths among nodes), and densities of distinctive patterns (e.g., total number of triangles).
We remark that, as different properties result in different admissible (feasible) network configurations, the NGP seeks configurations that are simultaneously feasible with respect to all specified properties.

In order to address the NGP, we will formulate optimization models according to two major aims.
The first is to concisely mathematically express each of the desired properties so that they can be conveniently evaluated given any network configuration that is explored during the solution process.
The second aim is to appropriately restrict the feasible solution space into only those network configurations that satisfy all property specifications.
Both aims are to be met with an appropriate selection of optimization variables and equality or inequality constraints, something that is discussed in detail in Sections~\ref{express} and~\ref{specify}.

Throughout the analysis, we preserve an overall MILP structure for our models and strive to involve expressions that have as tight LP-relaxations as possible.
We also focus on creating modular formulations, which we achieve by describing each property with an independent set of constraints that can be included as a building block to an overall formulation only when the
property in question is part of the problem specifications.
The modularity of our formulations is illustrated later in Section~\ref{studies}, where complete formulations are put together by utilizing the constraints presented as stand-alone examples in Sections~\ref{popular} and~\ref{enhancements}.

Let $V = \left\{1,2,\cdots,N\right\}$ be the set of the network's nodes, where $N$ is an {\em a priori} known constant - here the maximum allowable number of nodes.
Depending on the application, a set of properties may explicitly or implicitly impose that all $N$ nodes are non-isolated (i.e., have degree of at least one), or even that the network has a single component (i.e., is connected).
Let also $E=\left\{(i,j) \in V \times V : i < j\right\}$ be the set of all (undirected) edges that can possibly form between the nodes in $V$.
As the cases $N=0$ and $N=1$ are trivial for the NGP, we always assume for convenience that $N \geq 2$.

For ease of exposition, this paper focuses on undirected, simple networks.
Extensions to more general classes of networks can be handled by appropriately augmenting the
set $E$.
For example, loops can be considered by allowing $i =j$, while directed arcs can be introduced by
allowing $j > i$.
Furthermore, we consider nodes to be unlabeled;  that is, the nodes are permutable with
each other and the specification does not provide any particular reference to node labels.
We remark that the unlabeled case is in fact the more general and ``difficult'' variant of
the NGP, as achieving a more abstract specification must in general be more involved.

\subsection{Expressing Properties}\label{express}
The basic set of optimization variables utilized in our framework is a set of binary variables $x_{ij}$ to denote the existence of each edge $(i,j) \in E$:
\begin{equation}
x_{ij} = \left\{ \begin{array}{r@{\quad:\quad}l}
1 & \text{if edge $(i,j)$ exists in the solution}\\
0 & \text{otherwise.} \\
	     \end{array} \right.\label{xijdef}
\end{equation}

The $x$-variables suffice to fully define the configuration of the network.
Once their values are determined, every property of the network can be uniquely evaluated, and whether the property specifications are satisfied or not can be unambiguously determined.
However, only a handful of properties can be expressed in closed form as a linear function of these variables.
A notable example is the degree of a node, which can be expressed as the sum of all edges incident to that node.
Other properties will typically necessitate the introduction of additional auxiliary variables to facilitate their mathematical description.
For example, properties such as the global clustering coefficient require the encoding of whether particular features (in this case, triangles) have been formed or not, while expressing properties such as the network diameter would require the encoding of the shortest-path lengths among each pair of nodes.
To that end, it is often convenient to introduce a new set of variables $y_f$ to capture the formation of each feature $f \in F$, where $F$ is the set of features of interest, as well as a new set of variables $w_{ij}$ to express the length of the shortest path between each pair of nodes $(i,j) \in E$.
Restricting these variables to their intended values is discussed in Sections~\ref{features} and~\ref{shortestpaths}.

In general, most properties of interest can be expressed directly as a linear
combination of the $x$-, $y$- and $w$-variables:
\begin{equation}
p = \sum\limits_{(i,j) \in E} \alpha_{ij} x_{ij} + \sum\limits_{f \in F} \beta_f y_f + \sum\limits_{(i,j) \in E} \gamma_{ij} w_{ij} + \delta \label{pidef},
\end{equation}
where $\mathbf{\alpha}$, $\mathbf{\beta}$, $\mathbf{\gamma}$ and $\delta$ are parameters that depend on the type of property, and $p$ is an auxiliary variable to capture the property's value.

The use of an auxiliary $p$-variable helps the model from a human-readability and modularity perspective.
Furthermore, it facilitates the expression (in a linear fashion) of more complex, higher-level properties.
For example, once the local clustering coefficients of each node are represented by an appropriately defined set of $p$-variables, the average clustering coefficient of the entire network can be expressed based on their sum.
We note that, irrespectively of whether it represents a discrete or a continuous quantity, every $p$-variable may be declared to the MILP solver as a continuous variable within appropriate lower and upper bounds, $p^L$ and $p^U$.
These bounds can be readily derived based on interval arithmetic considerations.
We further remark that most commercial MILP solvers will automatically detect the special manner in which these additional variables participate in the formulation and remove them during a preprocessing phase.
Therefore, it is not expected that the computational performance of the solution process will be burdened by these types of modeling simplifications.
Examples of expressing a number of popular properties are provided in Section~\ref{popular}.

\subsubsection{Encoding the Existence of Features (Motifs)}\label{features}
The existence or absence in the designed network of characteristic patterns (e.g., triplets, triangles, cliques, stars) may be of interest as a result of some higher-level property (e.g., global clustering coefficient, clique number).
In such cases, the explicit encoding of the existence of these features is required to enable the expression of the desired property.
We observe that the formation of any feature (motif) corresponds to the concurrent existence and/or absence of certain edges.
For example, a triangle between nodes $1$, $2$ and $3$ exists if and only if all three edges $(1,2)$, $(1,3)$ and $(2,3)$ exist or, equivalently, if and only if $x_{12}x_{13}x_{23} = 1$.
Similarly, an open triplet along trajectory $2$--$1$--$3$ exists if and only if edges $(1,2)$ and $(1,3)$ exist, while edge $(2,3)$ does not, or if and only if $x_{12}x_{13}\left(1-x_{13}\right) = 1$.

For every feature $f \in F$, where $F$ is the set of all features of interest, let $E^{+}_f \subseteq E$ and $E^{-}_f \subseteq E$ be the corresponding collections of edges that need exist and be absent, respectively, and let a variable $y_f$ be defined as
\begin{equation}
y_{f} = \left\{ \begin{array}{r@{\quad:\quad}l}
1 & \text{if feature $f$ exists in the solution}\\
0 & \text{otherwise.} \\
	     \end{array} \right.\label{yfdef}
\end{equation}

In order to maintain linear structure in the model, the variable $y_f$ is not expressed directly as the product $\prod\limits_{(i,j) \in E^{+}_f} x_{ij} \prod\limits_{(i,j) \in E^{-}_f} \left(1-x_{ij}\right)$; rather it is restricted to an appropriate binary value through the following set of constraints:
\begin{eqnarray}
y_{f} &\leq& x_{ij} \;\;\;\;\;\;\;\;\;\;\;\;\; \forall (i,j) \in E^{+}_f \label{yf1a} \\
y_{f} &\leq& 1 - x_{ij} \;\;\;\;\;\; \forall (i,j) \in E^{-}_f \label{yf1b} \\
y_{f} &\geq& 1 - \sum\limits_{(i,j) \in E^{+}_f} \left(1-x_{ij}\right) - \sum\limits_{(i,j) \in E^{-}_f} x_{ij} \label{yf2} \\
y_{f} &\geq& 0 \label{yf3}.
\end{eqnarray}

We remark that, at any feasible solution (where all $x$-variables have attained binary values), the above constraints restrict the value of $y_f$ to either $0$ or $1$ without any explicit claim about its integrality.
This means that variables $y_f$ need not be defined as binary --rather as mere continuous-- in the context of the MILP solution algorithm.

An alternative to the $O (|E^{+}_f|+|E^{-}_f|)$ constraints~(\ref{yf1a}) and~(\ref{yf1b}) are the two aggregate constraints
\begin{eqnarray}
y_{f} &\leq& \frac{1}{|E^{+}_f|} \sum\limits_{(i,j) \in E^{+}_f} x_{ij} \label{yf1a'}\\
y_{f} &\leq& 1 - \frac{1}{|E^{-}_f|} \sum\limits_{(i,j) \in E^{-}_f} x_{ij}. \label{yf1b'}
\end{eqnarray}

Use of constraints~(\ref{yf1a})--(\ref{yf1b}) is preferred, as they result into tighter LP-relaxations.
Furthermore, their total-unimodularity endowes the overall LP-relaxations with an increased tendency to admit integral solutions.
The alternative constraints~(\ref{yf1a'})--(\ref{yf1b'}) are attractive in cases where the size of the formulation becomes an issue and has to be limited within an acceptable level.
We note that when using the aggregate constraints, binarity of variable $y_f$ is not guaranteed and must be explicitly declared to the MILP solver.
Examples of defining a number of common features are provided for illustration in Section~\ref{popular}.

\subsubsection{Encoding Shortest Paths}\label{shortestpaths}
A number of properties (e.g., betweeness-centrality, network diameter) depend on shortest-path information.
To that end, the explicit encoding of the shortest path between any node pair of interest $i \in V$ and $j \in V \setminus{\left\{i\right\}}$ is required.
Since a shortest path is itself an optimal quantity (out of all the paths that connect nodes $i$ and $j$, it is the one with the minimum length), it can be obtained via the solution of a mathematical optimization formulation.
Given a network configuration $\mathbf{x}$, the following formulation applies:
\begin{eqnarray}
&\min\limits_{\mathbf{f^{ij}}} & \sum\limits_{k \in V} \sum\limits_{l \in V} f^{ij}_{kl} \label{sp-p0} \\
& \text{s.t.} & f^{ij}_{kl} + f^{ij}_{lk} \leq x_{kl} \;\;\;\;\;\;\;\;\;\;\;\;\;\;\;\;\; \forall (k,l) \in E  \label{sp-p1} \\
&             & \sum\limits_{l \in V} f^{ij}_{kl} - \sum\limits_{l \in V} f^{ij}_{lk} = \xi^{ij}_k \;\;\;\;\;\; \forall k \in V \label{sp-p2} \\
&             & f^{ij}_{kk} = 0   \;\;\;\;\;\;\;\;\;\;\;\;\;\;\;\;\;\;\;\;\;\;\;\;\;\;\;\;\; \forall k \in V  \label{sp-p3} \\
&             & f^{ij}_{kl} \geq 0   \;\;\;\;\;\;\;\;\;\;\;\;\;\;\;\;\;\;\;\;\;\;\;\;\;\;\;\;\;\, \forall (k,l) \in V \times V \label{sp-p4},
\end{eqnarray}
where parameters $\xi^{ij}_i = +1$, $\xi^{ij}_j = -1$ and $\xi^{ij}_k = 0$, if $k \neq i,j$.

In the above formulation, each variable $f^{ij}_{kl}$ represents whether arc $(k \rightarrow l)$ participates or not in the shortest path between nodes $i$ and $j$.
Constraint~(\ref{sp-p1}) restricts participation to only those arcs that are eligible; that is, to only those arcs that correspond to edges that exist in the network.
The formulation ensures the integrality of the objective function at optimality.
To that end, it suffices that the $f$-variables be declared to the optimization solver as continuous.
We remark however that, due to the possibility of degeneracy (i.e., existence of more than one distinct shortest paths), some $f$-variables may admit fractional (i.e., neither $0$ nor $1$) values at optimality.
In cases when we are not only interested in the length of the shortest path but also interested in encoding the actual shortest path, the $f$-variables should be declared as binary.\footnote{For better performance of the branch-and-bound-based search process, it is advisable to deprioritize the $f$-variables during branching considerations.  By doing so, the $f$-variables will effectively be treated as continuous throughout the search process, and only at the very end (and only if necessary) will the solver take branching decisions on $f$ to resolve the degeneracy.}

An alternative formulation that can be used to obtain the shortest path between nodes $i$ and $j$ is the dual of the above (primal) formulation, which can be cast as follows:
\begin{eqnarray}
&\max\limits_{\mathbf{t^{ij}},\mathbf{u^{ij}},\mathbf{v^{ij}}} & \sum\limits_{(k,l) \in E} u^{ij}_{kl} + t^{ij}_i - t^{ij}_j \label{sp-d0} \\
& \text{s.t.} & \hspace{-0.23cm} \left.\begin{array}{rcl}
v^{ij}_{kl} + t^{ij}_{k} - t^{ij}_{l} &\leq& 1\\
v^{ij}_{kl} + t^{ij}_{l} - t^{ij}_{k} &\leq& 1\end{array}\right\} \;\;\;\; \forall (k,l) \in E  \label{sp-d1} \\
&             & u^{ij}_{kl} \leq v^{ij}_{kl} + U \left( 1 - x_{kl} \right) \;\;\;\;\;\; \forall (k,l) \in E  \label{sp-d2} \\
&             & t^{ij}_k \in {\mathbb R}   \;\;\;\;\;\;\;\;\;\;\;\;\;\;\;\;\;\;\;\;\;\;\;\;\;\;\;\;\;\;\; \forall k \in V  \label{sp-d3} \\
&             & u^{ij}_{kl}, v^{ij}_{kl} \leq 0   \;\;\;\;\;\;\;\;\;\;\;\;\;\;\;\;\;\;\;\;\;\;\;\;\; \forall (k,l) \in E \label{sp-d4},
\end{eqnarray}
where $\mathbf{t^{ij}}$, $\mathbf{u^{ij}}$ and $\mathbf{v^{ij}}$ are continuous variables, and $U=N-2$ is an appropriate ``big-M'' coefficient.

Although the aforementioned two formulations (the primal and the dual) are available to us for computing the shortest path between two given nodes of interest, none of the two can be readily utilized in the context of the NGP.
The complication arises due to the fact that the shortest path is only encoded when the formulations
are solved to optimality.
This cannot be guaranteed in general, because the objective function of the overall NGP formulation may have to be reserved for some other application-specific functional.
Furthermore, when shortest paths between additional pairs of nodes need to be encoded at the same time, as is often the case, then it would not be possible to concurrently optimize multiple interdependent --and generally competing-- objectives.
To address this shortfall, we use the strong duality theorem of linear programming, which ensures that when no gap exists between a feasible solution of the primal problem and a feasible solution of the dual problem (i.e., when these two solutions have equal objective values), then both solutions are the optimal solutions of their respective problems.
Given this fact, an explicit encoding of shortest-path information that does not utilize the overall formulation's objective function can be achieved by appending all constraints of both the primal and dual shortest path formulations (i.e., constraints~(\ref{sp-p1})--(\ref{sp-p4}) and~(\ref{sp-d1})--(\ref{sp-d4})), and in addition, appending the two equality constraints~(\ref{sp-o1}) and~(\ref{sp-o2}) to enforce that the objectives of the primal and dual formulations simultaneously attain a common value, $w_{ij}$.
\begin{eqnarray}
w_{ij} &=& \sum\limits_{k \in V} \sum\limits_{l \in V} f^{ij}_{kl} \label{sp-o1} \\
w_{ij} &=& \sum\limits_{(k,l) \in E} u^{ij}_{kl} + t^{ij}_i - t^{ij}_j \label{sp-o2},
\end{eqnarray}
where $w_{ij} \in \left[ 1, N-1 \right]$ are new variables that can be declared to the MILP solver as continuous variables (the constraints guarantee that they will attain discrete values).

We remark that an implication of encoding the shortest path between a pair of nodes $i$ and $j$ is the fact that the search is restricted to networks in which $i$ and $j$ participate in the same network component.\footnote{If shortest paths are encoded for all pairs $(i,j) \in E$, then the search is restricted to connected networks.}
This is due to the implicit requirement of the optimization process that any variable of the formulation attains a finite value.
%
However, in applications where the specified properties involve shortest paths (or some other property that depends on shortest paths), it is typically the case that the existence of the shortest paths is indeed desirable.

\subsection{Specifying Properties}\label{specify}
In the context of the NGP, the set of properties of interest should in general attain values within specified ranges.
Note that this generic specification admits also the cases where there is only a single (either lower or upper) required bound, as well as the case where a property must equal some specified value exactly.
In the following, we describe how range restrictions can be imposed for linear, fractional, and sequence properties, as well as restrictions on collective statistics of those.
Examples of imposing restrictions for a number of popular properties will be presented in Section~\ref{slacks}.
We remark that the larger the number of property specifications, the smaller the feasible solution space is.
To that end, property combinations should be chosen judiciously so as not to inadvertedly correspond to infeasible problem instances (i.e., instances where no feasible network configuration exists).

\subsubsection{Linear Properties}\label{specify_lp}
We define a linear property $lp$ to be any property that can be directly expressed as a linear combination of any of the formulation's declared variables (e.g., the total number of triplets or triangles in the network can be expressed as sums of appropriately defined $y$-variables).
This definition obviously admits any property that can be expressed as a $p$-variable according to constraint~(\ref{pidef}).
Specifying that $lp$ attains some value in the range $\left[\pi^{L},\pi^{U}\right]$ can be trivially achieved by appending in the overall formulation the explicit bounds
\begin{equation}
\pi^L \leq lp \leq \pi^U \label{lp_bounds}.
\end{equation}
It is assumed that the linear equality constraint that constitutes the definition of property $lp$ has also been included in the formulation.

\subsubsection{Fractional Properties}\label{specify_fp}
We define a fractional property $fp$ to be any property that can be expressed as the ratio of two linear properties, $fp = \frac{lp^A}{lp^B}$, where the linear property $lp^B$ cannot attain the value of zero.
An example of a fractional property would be the global clustering coefficient, which is defined as thrice the ratio between the total triangles and triplets in the network.
Specifying that a fractional property attains some value in the range $\left[\pi^{L},\pi^{U}\right]$ can be achieved by appending in the overall formulation the constraints\footnote{Constraints~(\ref{fp_bounds}) apply for strictly positive $lp^B$.
If instead $lp^B$ is strictly negative, the inequality signs should be reversed.}
\begin{equation}
\pi^L lp^B \leq lp^A \leq \pi^U lp^B \label{fp_bounds}.
\end{equation}
It is assumed that the linear equality constraints that constitute the definitions of properties $lp^A$ and $lp^B$ have also been included in the formulation.

\subsubsection{Sequence Properties}\label{specify_seq}
Most specified properties are scalar quantities (e.g., the clustering coefficient of a node, the girth of the network).
However, a property may also be a vector; a case in point would be the scalar values that a collection of nodes attains for some other, scalar property.
For example, the degree sequence of all neighbors of a given node $i$ is a vector property.
If we want to impose it on the network to be constructed, we need to have explicitly encoded which node $j$ (that is a neighbor of node $i$) attains which value from the sequence.

More formally, let a node-specific property of interest and let variables $p_i$ be this property's values for each node $i$, given a network configuration $\mathbf{x}$.
It is assumed that the $p$-variables attain appropriate values at any feasible solution (e.g., as a result of constraints of type~(\ref{pidef}) having been included in the formulation).
In addition, let a subset of $M$ nodes ($M \leq N$), and let participation of a node $i \in V$ in this subset be dictated by an appropriately defined binary variable $z_i$; that is, $z_i = 1$, if node $i$ participates in the subset, and $z_i = 0$, if it does not, such that $\sum\limits_{i \in V} z_i = M$ (the activation of the $z$-variables is formally discussed in Section~\ref{subsets}).
Finally, let a sequence of ranges $\mathbf{\Pi} = \left\{ \left[\pi_1^L,\pi_1^U\right], \left[\pi_2^L,\pi_2^U\right], \cdots, \left[\pi_M^L,\pi_M^U\right] \right\}$ that is to be specified on the node subset for the property of interest.
In other words, it is desired that the solution exhibits a one-to-one correspondence between the $p_i$ values for the $M$ nodes and the ranges in the sequence $\mathbf{\Pi}$.
We note that the admissible correspondence may not be unique, if any two ranges in $\mathbf{\Pi}$ overlap.
We introduce a new set of binary variables $q_{im}$ to denote the assignment of nodes to elements of the sequence:
\begin{equation}
q_{im} = \left\{ \begin{array}{r@{\quad:\quad}l}
1 & \text{if node $i$ is assigned to element $m$}\\
0 & \text{otherwise,} \\
	     \end{array} \right.\label{qimdef}
\end{equation}
where $z_i = 1$ and $\pi_m^L \leq p_i \leq \pi_m^U$ are necessary conditions for assignment.

The definition of the $q$-variables and the specification of the sequence is enforced via the following set of constraints:
\begin{eqnarray}
&\left.\begin{array}{rcl}
p_i &\geq& p_i^L - \left(p_i^L - \pi_m^L \right) q_{im} \\
p_i &\leq& p_i^U - \left(p_i^U - \pi_m^U \right) q_{im}\end{array}\right\}& \forall (i,m) \in V \times \left\{1,2,\cdots,M\right\} \label{qim1} \\
& \sum\limits_{i \in V} q_{im} =  1& \forall m \in \left\{1,2,\cdots,M\right\} \label{qim2} \\
& \sum\limits_{m = 1}^M q_{im} \leq z_{i}&  \forall i \in V \label{qim3}\\
& \sum\limits_{m = 1}^M q_{im} \geq z_{i}&  \forall i \in V \label{qim4},
\end{eqnarray}
where $p_i^L$ and $p_i^U$ are appropriate lower and upper limits for the value of variable $p_i$.

We remark that in cases where the sequence refers to nodes belonging to the intersection of multiple subsets, the $q$-variables need to be appropriately restricted.
This can be achieved by (i) appending separate constraints~(\ref{qim3}) for each applicable subset and (ii) replacing the right-hand side in constraints~(\ref{qim4}) with an aggregate term $1 - \sum\limits_e  \left(1 - z_i^e\right)$ ($e$ is used here to index over the collection of subsets of interest and $z_i^e$ are the corresponding $z$-variables for subset $e$).

\subsubsection{Statistical Properties}\label{specify_st}
Instead of specific values or sequences, one may want to specify collective statistics,
such as the arithmetic mean or the median of a given property.
Typically these statistical metrics are to be evaluated across the complete set of nodes (e.g., the case of the network's average clustering coefficient).
%
However, in certain contexts they may refer only to a subset of nodes, which could itself be a variable quantity that depends on the network configuration $\mathbf{x}$ (e.g., the case of the average degree among the neighbors of a given node).

More formally, let a node-specific property of interest and let variables $p_i$ attain appropriate property values for each node $i \in V$.
In addition, let variables $z_i$ define participation in a node subset of interest.
To assist us in the expression of statistical metrics, we introduce a new set of continuous variables $p'_i$ to be equal to $p_i z_i + \tilde{p} \left(1 - z_i \right)$; that is, to be equal to the value of property $p_i$, if $z_i = 1$, or equal to some constant value $\tilde{p}$, if $z_i = 0$.
This definition of the $p'$-variables is enforced with the following set of constraints:
\begin{eqnarray}
&\left( p_i^L - \tilde{p} \right) z_i \leq p'_i - \tilde{p} \leq \left( p_i^U - \tilde{p} \right) z_i& \;\;\;\; \forall i \in V \label{pip1}\\
&\left( p_i^L - \tilde{p} \right) \left(1 - z_i\right) \leq p_i - p'_i \leq \left( p_i^U - \tilde{p} \right) \left(1 - z_i\right) & \;\;\;\; \forall i \in V \label{pip2},
\end{eqnarray}
where $p_i^L$ and $p_i^U$ are appropriate lower and upper limits for the value of variable $p_i$.  We remark that, as a result of these constraints, $p'_i \in \left[ \min\left\{\tilde{p}, p_i^L\right\}, \max\left\{\tilde{p}, p_i^U\right\}\right]$, for all $i \in V$.

Furthermore, we note that in cases where the sequence refers to nodes belonging to the intersection of multiple subsets, the $p'$-variables need to be appropriately restricted.
This can be achieved by (i) appending separate constraints~(\ref{pip1}) for each applicable subset and (ii) replacing the term $\left(1 - z_i\right)$ in constraints~(\ref{pip2}) with an aggregate term $\sum\limits_e  \left(1 - z_i^e\right)$ ($e$ is used here to index over the collection of subsets of interest and $z_i^e$ are the corresponding $z$-variables for subset $e$).

Given the constraint set~(\ref{pip1}) and~(\ref{pip2}), the sum of the property values across the subset of interest can be linearly expressed by setting $\tilde{p} = 0$ and summing $p'$-variables across the constant node set $V$, $\sum\limits_{i \in V} p'_i$.
The mean across the subset of interest can then be expressed as the ratio $\sum\limits_{i \in V} p'_i / \sum\limits_{i \in V} z_i$.
Therefore, the mean constitutes a fractional property and specifying it can be achieved according to the discussion in Section~\ref{specify_fp}.
Higher-order moments could be specified along the same lines, assuming integral powers of the properties ($p_i^2$, $p_i^3$, etc.) can also be expressed in a linear fashion.
We remark that the specification of higher-order moments simplifies if the moments are taken around a constant value, or if the corresponding lower-order moments are also specified to fixed values.

We will conclude this discussion by addressing the case of specifying the
median\footnote{Given an ordered set of $n$ values $\alpha_i$, $i \in \left\{1,2,\cdots,n\right\}$,
we define the median to be the value $\alpha_{\left(n+1\right)/2}$, if $n$ is odd,
or any value in the range $\left[\alpha_{n/2},\alpha_{\left(n+2\right)/2}\right]$, if $n$ is even.}
$\hat{p}$ of a set of property values $p_i$, $i \in V \setminus \left\{i : z_i = 0\right\}$.
We will show that the median can be
linearly expressed in the formulation and, thus, can be directly specified as a linear property.
We first note that the values $\hat{p}^L = \min\limits_{i \in V} \left\{p_i^L\right\}$ and $\hat{p}^U = \max\limits_{i \in V} \left\{p_i^U\right\}$ are appropriate lower and upper limits for the value of the median, irrespectively of the actual participation of nodes in the subset of interest.
We then define $p'$-variables based on $\tilde{p} = \hat{p}^U$ and introduce two new sets of binary variables, $\mathbf{r^{+}}$ and $\mathbf{r^{-}}$, to indicate respectively whether the median is larger or smaller than each of the $p'$-variables.
Similarly to~(\ref{zi}), the following constraints enforce the definitions of the $r$-variables:
\begin{eqnarray}
&\left( p_i^L - \hat{p}^U \right) r_i^{+} \leq p'_i - \hat{p} \leq \left( p_i^U - \hat{p}^L \right) \left(1 - r_i^{+}\right)& \label{rip}\\
&\left( p_i^U - \hat{p}^L \right) r_i^{-} \geq p'_i - \hat{p} \geq \left( p_i^L - \hat{p}^U \right) \left(1 - r_i^{-}\right)& \label{rim}.
\end{eqnarray}

Given this definition, the median $\hat{p}$ will attain the proper value as a result of including in the formulation the inequalities
\begin{eqnarray}
&\sum\limits_{i \in V} z_i \leq 2 \sum\limits_{i \in V} r_i^{+} \leq 1 + \sum\limits_{i \in V} z_i& \label{medianp}\\
&2N - \sum\limits_{i \in V} z_i \leq 2 \sum\limits_{i \in V} r_i^{-} \leq 1 + 2N - \sum\limits_{i \in V} z_i& \label{medianm}\\
&\sum\limits_{i \in V} z_i \geq 1 ,& \label{median1}
\end{eqnarray}
where the last constraint requires explicitly that the subset of interest is not empty.  We note that the LP-relaxation can be tightened by adding constraints $r_i^{-}~\geq~1-z_i$ and $r_i^{+}~\leq~z_i$ for all $i \in V$.

\subsubsection{Participation of Nodes in Subsets}\label{subsets}
Every property refers to an appropriately defined subset of the network's nodes.
Most often properties are node-specific; that is, given any network configuration $\mathbf{x}$, each node $i \in V$ attains a separate value for this property (e.g., the properties of degree and clustering coefficient).
In this case, the applicable subset consists of a single node.
Properties may also refer to the network as a whole (e.g., the degree sequence), in which case the subset consists of the complete set $V$.
There are additional cases, however, where a property may refer to some other subset of the nodes (e.g., the case of specifying the maximum shortest-paths among the nodes with a given clustering coefficient).
We note that, in the general case, the participation of nodes in a specified subset of interest is itself variable; that is, it depends on the network configuration $\mathbf{x}$ and it is not a-priori known which nodes belong to the subset.
In such cases, whether a node belongs or not to the subset needs to be appropriately encoded in the formulation.

Given a node-specific property, we define a subset of interest as the collection of those nodes $i \in V$ for which the corresponding values of this property, $p_i$, evaluate above some specified threshold $\pi$.
To indicate whether each node $i \in V$ belongs or not to the subset of interest, we introduce a new set of binary variables
\begin{equation}
z_{i} = \left\{ \begin{array}{r@{\quad:\quad}l}
1 & \text{if node $i$ participates in the subset}\\
0 & \text{otherwise,} \\
	     \end{array} \right.\label{zidef}
\end{equation}
where $p_i \geq \pi$ is the necessary and sufficient condition for participation.

The definition of the $z$-variables is enforced via the constraints\footnote{We note that constraints~(\ref{zi}) do not completely enforce the definition~(\ref{zidef}), as variable $z_i$ remains unrestricted in the case when variable $p_i$ exactly attains the value $\pi$.  However, in cases where $p_i$ attains values from a discrete set, one can always identify some small $\epsilon$ to subtract from $\pi$ in~(\ref{zi}) so that the equality $p_i = \pi$ induces the restriction $z_i = 1$.  In cases where $p_i$ attains values from a continuum, there is no practical concern, as there is zero likelihood that $p_i$ attains any specific value exactly.}
\begin{equation}
\left( p_i^L - \pi \right) \left(1 - z_i\right) \leq p_i - \pi \leq \left( p_i^U - \pi \right) z_i \label{zi},
\end{equation}
where $p_i^L$ and $p_i^U$ are appropriate lower and upper limits for the value of variable $p_i$.
It is assumed that the $p$-
variables attain appropriate values as a result of constraints~(\ref{pidef}) having been included in the formulation.

We note that the subset definition presented above can readily admit the case where $\pi$ refers to an upper limit; that is, the case where participation in the subset is realized when $p_i \leq \pi$.
In such a case, one need use $-p_i$ as the property and $-\pi$ as the threshold value, which is equivalent to performing the substitution $z_i \leftarrow \left(1 - z_i \right)$ in constraints~(\ref{zi}).
More complex subsets can then be expressed as intersections of these two simpler cases.
For example, a subset that is defined based on a property evaluating within some specified range, $\pi^L \leq p_i \leq \pi^U$, can be described as the intersection of the two simple subsets $p_i \geq \pi^L$ and $p_i \leq \pi^U$.

\section{Popular Properties}\label{popular}
In this section, we use the modeling elements presented in Section~\ref{modeling} in order to describe a number of popular network properties and associated features, which can then be used in various specifications for the NGP.
Although not an exhaustive list of properties, the compilation of examples illustrates the generality of the proposed MILP-based approach.
We begin by illustrating how a number of common features can be defined.
\begin{example}
{\bf (2-Path).} Let variable $ytp_{ijk}$ indicate whether a 2-path (i.e., open or closed triplet) exists along trajectory $j$--$i$--$k$, where $i \in V$ and $(j,k) \in \left\{\left(V \setminus{\left\{i\right\}}\right)^2 : j < k\right\}$.  This variable can be defined with the constraint set
\begin{equation}
\begin{array}{c}
ytp_{ijk} \leq x_{j'j''} \;;\;ytp_{ijk} \leq x_{k'k''} \\
ytp_{ijk} \geq x_{j'j''} + x_{k'k''} - 1\\
ytp_{ijk} \geq 0,
\end{array}\label{ytp}
\end{equation}
where $j' = \min\left\{i,j\right\}$, $j'' = \max\left\{i,j\right\}$, $k' = \min\left\{i,k\right\}$ and $k'' = \max\left\{i,k\right\}$.
Generalizations to longer paths can be readily addressed. \qed
\end{example}
\begin{example}
{\bf (3-Cycle).} Let variable $ytr_{ijk}$ indicate whether a 3-cycle (i.e., triangle) exists between three nodes $(i,j,k) \in T$, where $T = \left\{(i,j,k) \in V^3 : i < j < k\right\}$.  This variable can be defined with the constraint set
\begin{equation}
\begin{array}{c}
ytr_{ijk} \leq x_{ij} \;;\;ytr_{ijk} \leq x_{ik} \;;\;ytr_{ijk} \leq x_{jk} \\
ytr_{ijk} \geq x_{ij} + x_{ik} + x_{jk} - 2\\
ytr_{ijk} \geq 0.
\end{array}\label{ytr}
\end{equation}
Generalizations to longer cycles (i.e., squares) can be readily addressed. \qed
\end{example}
\begin{example}
{\bf (4-Clique).} Let variable $yclq_{ijkl}$ indicate whether a clique exists between four nodes $(i,j,k,l) \in \left\{V^4 : i < j < k < l\right\}$.  This variable can be defined with the constraint set
\begin{equation}
\begin{array}{c}
yclq_{ijkl} \leq x_{ij} \;;\;yclq_{ijkl} \leq x_{ik} \;;\;yclq_{ijkl} \leq x_{il} \\
yclq_{ijkl} \leq x_{jk} \;;\;yclq_{ijkl} \leq x_{jl} \;;\;yclq_{ijkl} \leq x_{kl} \\
yclq_{ijkl} \geq x_{ij} + x_{ik} + x_{il} + x_{jk} + x_{jl} + x_{kl} - 5\\
yclq_{ijkl} \geq 0,
\end{array}\label{yclq}
\end{equation}
or equivalently, with the constraint set
\begin{equation}
\begin{array}{c}
yclq_{ijkl} \leq \frac{1}{6} \left(x_{ij} + x_{ik} + x_{il} + x_{jk} + x_{jl} + x_{kl}\right)\\
yclq_{ijkl} \geq x_{ij} + x_{ik} + x_{il} + x_{jk} + x_{jl} + x_{kl} - 5\\
yclq_{ijkl} \in \left\{0,1\right\}.
\end{array}\label{yclq'}
\end{equation}
Generalizations to larger cliques can be readily addressed. \qed
\end{example}
\begin{example}
{\bf (4-Star).} Let variable $ystr_{ijkl}$ indicate whether node $i \in V$ corresponds to the center of a star involving nodes $(j,k,l) \in \left\{\left(V \setminus{\left\{i\right\}}\right)^3 : j < k < l\right\}$.  This variable can be defined with the constraint set
\begin{equation}
\begin{array}{c}
ystr_{ijkl} \leq x_{j'j''} \;;\;ystr_{ijkl} \leq x_{k'k''} \;;\;ystr_{ijkl} \leq x_{l'l''} \\
ystr_{ijkl} \leq 1 - x_{jk} \;;\;ystr_{ijkl} \leq 1 - x_{jl} \;;\;ystr_{ijkl} \leq 1 - x_{kl} \\
ystr_{ijkl} \geq x_{j'j''} + x_{k'k''} + x_{l'l''} - x_{jk} - x_{jl} - x_{kl} - 2\\
ystr_{ijkl} \geq 0,
\end{array}\label{ystr}
\end{equation}
or equivalently, with the constraint set
\begin{equation}
\begin{array}{c}
ystr_{ijkl} \leq \frac{1}{3} \left(x_{j'j''} + x_{k'k''} + x_{l'l''} \right) \\
ystr_{ijkl} \leq 1 - \frac{1}{3} \left(x_{jk} + x_{jl} + x_{kl} \right) \\
ystr_{ijkl} \geq x_{j'j''} + x_{k'k''} + x_{l'l''} - x_{jk} - x_{jl} - x_{kl} - 2\\
ystr_{ijkl} \in \left\{0,1\right\},
\end{array}\label{ystr'}
\end{equation}
where $j' = \min\left\{i,j\right\}$, $j'' = \max\left\{i,j\right\}$, $k' = \min\left\{i,k\right\}$, $k'' = \max\left\{i,k\right\}$, $l' = \min\left\{i,l\right\}$ and $l'' = \max\left\{i,l\right\}$.
Generalizations involving more nodes can be readily addressed. \qed
\end{example}

Given such definitions of features, a number of properties can be expressed in a linear or fractional fashion.  Some examples follow.
\begin{example}
{\bf (Degree).} The degree of each node $i \in V$ can be captured by the variable
\begin{equation}
pd_i = \sum\limits_{(j,k) \in E: \atop \left\{(j = i) \vee (k = i)\right\}} x_{jk}.  \label{pd}
\end{equation} \qed
\end{example}
\begin{example}
{\bf (Number of Triplets).} The number of open and closed triplets (i.e., 2-paths) centered at each node $i \in V$ can be captured by the variable
\begin{equation}
pntp_i = \sum\limits_{(j,k) \in E} ytp_{ijk}.  \label{pntp} \qed
\end{equation} 
\end{example}
\begin{example}
{\bf (Number of Triangles).} The number of triangles (i.e., 3-cycles) in which each node $i \in V$ participates can be captured by the variable
\begin{equation}
pntr_i = \sum\limits_{(j,k,l) \in T: \atop \left\{(j = i) \vee (k = i) \vee (l = i)\right\}} ytr_{jkl}. \label{pntr} \qed
\end{equation} 
\end{example}
\begin{example}
{\bf (Clustering Coefficient).} The (local) clustering coefficient of a node $i \in V$ can be captured by the variable
\begin{equation}
pcc_i = \frac{pntr_i}{pntp_i} \label{pcc}
\end{equation}
and, thus, constitutes a fractional property.  However, when node $i$ has a specified (fixed) degree $d_i$, then the local clustering coefficient can also be captured by the linear expression
\begin{equation}
pcc_i = \frac{1}{\left(d_i \atop 2 \right)} pntr_i, \label{pcc_fix}
\end{equation}
where $\left(d_i \atop 2 \right)$ is the constant number of 2-paths centered at node $i$. \qed
\end{example}
\begin{example}
{\bf (Average Clustering Coefficient).} The average clustering coefficient can be captured by the variable
\begin{equation}
pacc = \frac{1}{N} \sum\limits_{i \in V} pcc_i. \label{pacc} \qed
\end{equation} 
\end{example}
\begin{example}
{\bf (Global Clustering Coefficient).} The global clustering coefficient of a network can be captured by the variable
\begin{equation}
pgcc = \frac{\sum\limits_{i \in V} pntr_i}{\sum\limits_{i \in V} pntp_i} \label{pgcc}
\end{equation}
and, thus, constitutes a fractional property.
However, when the network has a specified (fixed) degree sequence $\mathbf{d} = \left\{d_1, d_2, \cdots, d_N\right\}$, then the global clustering coefficient can also be captured by the linear expression
\begin{equation}
pgcc = \frac{3}{\sum\limits_{i \in V} \left(d_i \atop 2 \right)} \sum\limits_{i \in V} pntr_i, \label{pgcc_fix}
\end{equation}
where $\sum\limits_{i \in V} \left(d_i \atop 2 \right)$ is the constant number of total 2-paths in the network. \qed
\end{example}
\begin{example}
{\bf (Closeness Centrality).} The inverse of the closeness centrality of each node $i \in V$ can be captured by the variable
\begin{equation}
piclc_i = \frac{1}{N-1} \sum\limits_{(j,k) \in E: \atop \left\{(j = i) \vee (k = i)\right\}} w_{jk}.  \label{piclc}
\end{equation}
Specifying bounds $\left[\pi^L,\pi^U\right]$ for the closeness centrality of a node can be achieved through specifying bounds $\left[\frac{1}{\pi^U},\frac{1}{\pi^L}\right]$ for its inverse. \qed
\end{example}
\begin{example}
{\bf (Average Path Length).} The average path length of a network can be captured by the variable
\begin{equation}
papl = \frac{1}{|E|} \sum\limits_{(i,j) \in E} w_{ij}, \label{papl}
\end{equation}
where $|E| = N \left(N-1\right)/2$ is the cardinality of the set of edges $E$. \qed
\end{example}
\begin{example}
{\bf (Characteristic Path Length).} The characteristic path length of a network (i.e., defined here as the median of all its shortest paths) can be captured by a variable $pcpl$, which assumes its appropriate value via the constraints\footnote{Constraints~(\ref{pcpl}) are the equivalent of constraints~(\ref{medianp}) and~(\ref{medianm}) for the case of a fixed subset of interest with cardinality $|E|$ (in particular, the set of edges $E$).}
\begin{equation}
\begin{array}{c}
\hspace{-0.23cm} \left.\begin{array}{c}
\begin{array}{rcr}
-U_1 \; rcpl_{ij}^{+} &\leq w_{ij} - pcpl \leq&  U_2 \left(1 - rcpl_{ij}^{+}\right)\\
 U_2 \; rcpl_{ij}^{-} &\geq w_{ij} - pcpl \geq& -U_1 \left(1 - rcpl_{ij}^{-}\right)
\end{array}\\
rcpl_{ij}^{+},rcpl_{ij}^{-} \in \left\{0,1\right\}
\end{array}\right\} \;\;\; \forall (i,j) \in E \\
2\sum\limits_{(i,j) \in E} rcpl_{ij}^{+} = 2\sum\limits_{(i,j) \in E} rcpl_{ij}^{-} = |E| + \mathbbm{1}_{\left\{|E| \text{ is odd}\right\} },
\end{array} \label{pcpl}
\end{equation}
where $U_1 = \left \lceil{ \frac{N}{2} }\right \rceil - 1$ and $U_2 = N - 2$ are appropriate ``big-M'' coefficients, $|E| = N \left(N-1\right)/2$ is the cardinality of the set of edges $E$, and $\mathbbm{1}_{\left\{\cdot\right\}}$ is the indicator function. \qed
\end{example}
\begin{example}
{\bf (Sum of Degrees of Neighbors).} The sum of degrees of all neighbors of each node $i \in V$ can be captured by variables $psdn_i$, which assume their appropriate values via the constraints
\begin{equation}
\begin{array}{cl}
\left.\begin{array}{rcl}
d^L \; x_{kl} \leq &pd'_{kl} &\leq d^U \; x_{kl}\\
d^L \; \left( 1- x_{kl}\right) \leq &pd_l - pd'_{kl} &\leq d^U \; \left(1 - x_{kl}\right)\\
d^L \; x_{kl} \leq &pd'_{lk} &\leq d^U \; x_{kl}\\
d^L \; \left( 1- x_{kl}\right) \leq &pd_k - pd'_{lk} &\leq d^U \; \left(1 - x_{kl}\right)
\end{array}\right\} & \forall (k,l) \in E \\
\hspace{-0.70cm} psdn_i = \sum\limits_{j \in V} pd'_{ij} & \forall i \in V,
\end{array} \label{psdn}
\end{equation}
where $pd'_{kl}$, $pd'_{lk}$ are $p'$-variables defined for each $(k,l) \in E$ with $\tilde{p}=0$ and with $x_{kl}$ serving as the corresponding $z$-variables, and $d^L$, $d^U$ are applicable lower and upper bounds for the degrees of each node $i \in V$ (such that $0 \leq d^L \leq d^U \leq N-1$).
However, when the network has a specified (fixed) degree sequence $\mathbf{d} = \left\{d_1, d_2, \cdots, d_N\right\}$, one can instead use the simpler expression
\begin{equation}
psdn_i = \sum\limits_{j \in V : \atop \left\{(i,j) \in E\right\}} d_j x_{ij} + \sum\limits_{j \in V : \atop \left\{(j,i) \in E\right\}} d_j x_{ji}.\label{psdn_fix} \qed 
\end{equation} 
\end{example}
\begin{example}
{\bf (Number of Nodes with Given Degree).} Given applicable lower and upper bounds $d^L$ and $d^U$ (such that $0 \leq d^L \leq d^U \leq N-1$) for the degree of each node in the network, the number of nodes that attain each possible degree $q \in Q$, where $Q = \{d^L,d^L+1,\cdots,d^U-1,d^U\}$, can be captured by variables $pnnd_q$.
These variables assume their appropriate values via the constraints
\begin{equation}
\left.\begin{array}{c}
\left.\begin{array}{c}
\left( d^L \; - q + \epsilon \right) \left(1 - z_{qi}\right) \leq pd_i - q + \epsilon \leq \left( d^U - q + \epsilon\right) z_{qi}\\
z_{qi} \in \left\{0,1\right\}
\end{array} \right\} \; \forall i \in V \\
\hspace{-1.75cm} pnnd_q = \sum\limits_{i \in V} z_{qi} - \sum\limits_{i \in V} z_{\left(q+1\right)i}\\
\end{array} \right\} \; \forall q \in Q, \label{pnnd}
\end{equation}
where $z_{qi}$, $q \in Q$, are $z$-variables that are activated when $pd_i \geq q$ and deactivated when $pd_i < q$,
$z_{\left(d^U+1\right)i} = 0$ for all $i \in V$ are parameters defined for notational convenience,
and $\epsilon \in \left(0,1\right)$ is a small parameter to induce activation of variable $z_{qi}$ when $pd_i$ equals $q$ exactly. \qed
\end{example}

\section{Algorithmic Enhancements}\label{enhancements}
As will become evident with the case studies we will present in Section~\ref{studies}, the modeling techniques discussed so far provide considerable amount of flexibility in putting together NGP formulations for a wide variety of property specifications of interest.
However, the NGP is in general a highly-complex combinatorial problem and solving it via standard optimization algorithms may be challenging for all but the smallest instances.
In this section, we discuss a number of algorithmic enhancements that we have found to improve the performance of solution algorithms in the NGP setting.
More specifically, in Section~\ref{slacks} we discuss the utilization of an appropriate objective function that attempts to minimize the cumulative deviation of the solution from the allowable ranges dictated by the specifications.
In Section~\ref{symmbreak}, we discuss the issue of breaking the symmetry embedded inherently in any NGP formulation by appropriately restricting the search space to solutions that are non-isomorphic to each other.

\subsection{Deviation of Properties from Specified Values}\label{slacks}
The discussion in Section~\ref{modeling} addressed efficient ways to express and specify (restrict) properties of interest in an overall MILP formulation.
In all these cases, setting up the modeling of the problem did not require
the formulation's objective function.
This allows us to formulate application-specific objective functions.
For example, out of all those network configurations $\mathbf{x}$ that meet a certain specification, one may seek to identify the network that maximizes or minimizes some particular property, in which case the objective must be the appropriate expression quantifying that property.

Nevertheless, there are many applications for which an explicit metric to assess the fitness of a network configuration does not exist, and the problem at hand is merely a feasibility problem; that is, to identify one or more configurations $\mathbf{x}$ that meet the specification, without consideration to how ``good'' or ``bad'' each of the obtained solutions are.
In such cases, a suitable objective may be devised to enhance the computational performance of the solution process.
In particular, variables can be introduced to encode the deviation of each of the properties from its nearest acceptable (as per the specification) value.
The objective would then require that the sum of all these ``slack'' variables be minimized.
Such an objective quantifies the total deviation from the feasible space of the problem and allows the solver to efficiently guide the search and locate feasible solutions.

As discussed, the property specifications come in the general form $g_m^L \leq g_m(\mathbf{x}) \leq g_m^U$, where $m$ is used here to index over each specified property, $g_m$ is the corresponding functional of the formulation's variables and $g_m^L$, $g_m^U$ are appropriate constants (see, e.g., the constraints~(\ref{lp_bounds}), (\ref{fp_bounds}), (\ref{qim1}), and~(\ref{symbreaklite1})).
Continuous slack variables $\mathbf{s}$ are introduced in the formulation in a way that relaxes each of these specification constraints:
\begin{eqnarray}
&g_m^L \leq g_m(\mathbf{x}) + s_m^{-} - s_m^{+}\leq g_m^U& \label{s1}\\
&s_m^{-}, s_m^{+} \geq 0,& \label{s2}
\end{eqnarray}
where $s_m^{-}$ and $s_m^{+}$ are pairs of variables to encode respectively the negative and positive deviation of quantities $g_m$ from their allowable ranges $\left[g_m^L,g_m^U\right]$.
The sum of all $s$-variables, $\sum\limits_m \left(s_m^{-} + s_m^{+}\right)$, should constitute the objective function to be minimized.

We remark that for a feasible solution all slack variables need to be zero at optimality.\footnote{Note that corresponding variables $s^{-}$ and $s^{+}$ cannot be both non-zero at an optimal solution.}
An optimal solution with a strictly positive objective (as defined above) corresponds to an infeasible NGP, and would signify that the property specification is not attainable.
Examples of specifying popular properties via use of $s$-variables are presented below.
Whenever the corresponding specifications are required, these constraints are to be appended during the modular build-up of a complete formulation.
In Section~\ref{studies}, we present a number of illustrative formulations, and we elaborate further on how these constraints contribute in the context of a complete optimization model.


\begin{example}
{\bf (Specification of Degree).} Let $pd_i$ be the variable capturing the degree of a node $i \in V$ as per constraint~(\ref{pd}).
A lower bound $d^L$ and an upper bound $d^U$ can be specified on the degrees of this node via the constraints
\begin{equation}
d^L \leq pd_i + sd_i^{-} - sd_i^{+} \leq d^U \label{spec-deg},
\end{equation}
where $sd_i^{-}, sd_i^{+} \geq 0$ are a pair of slack variables to be minimized.\qed
\end{example}
\begin{example}
{\bf (Specification of Degree Sequence).} Let $pd_i$ be the variables capturing the degree of each node $i \in V$ as per constraint~(\ref{pd}).
A degree sequence $\mathbf{d} = \left\{d_1, d_2, \cdots, d_N\right\}$ can be specified via the constraints
\begin{equation}
d_i \leq pd_i + sd_i^{-} - sd_i^{+} \leq d_i \;\; \forall i \in V \label{spec-degseq},
\end{equation}
where $sd_i^{-}, sd_i^{+} \geq 0$ are two sets of slack variables to be minimized.\qed
\end{example}
\begin{example}
{\bf (Specification of Average Clustering Coefficient).} Let $pacc$ be the variable capturing the average clustering coefficient as per constraint~(\ref{pacc}).
A lower bound $acc^L$ and an upper bound $acc^U$ can be specified on the average clustering coefficient via the constraints
\begin{equation}
acc^L \leq pacc + sacc^{-} - sacc^{+} \leq acc^U \label{spec-acc},
\end{equation}
where $sacc^{-}$, $sacc^{+} \geq 0$ are a pair of slack variables to be minimized.\qed
\end{example}
\begin{example}
{\bf (Specification of Global Clustering Coefficient).}  Let the global clustering coefficient of a network as defined in constraint~(\ref{pgcc}).
Since this quantity constitutes a fractional property, a lower bound $gcc^L$ and an upper bound $gcc^U$ can be specified via the constraints
\begin{equation}
gcc^L \sum\limits_{i \in V} pntp_i \leq \sum\limits_{i \in V} pntr_i + sgcc^{-} - sgcc^{+} \leq gcc^U \sum\limits_{i \in V} pntp_i \label{spec-gcc},
\end{equation}
where $sgcc^{-}$, $sgcc^{+} \geq 0$ are a pair of slack variables to be minimized.

However, when the network has a specified (fixed) degree sequence, then the global clustering coefficient can be specified as a linear property.
Let $pgcc$ be the variable capturing the global clustering coefficient as per constraint~(\ref{pgcc_fix}).  We can directly utilize this variable in the constraints
\begin{equation}
gcc^L \leq pgcc + sgcc^{-} - sgcc^{+} \leq gcc^U \label{spec-gcc_fix},
\end{equation}
where $sgcc^{-}$, $sgcc^{+} \geq 0$ are a pair of slack variables to be minimized.\qed
\end{example}
\begin{example}
{\bf (Specification of Average Path Length).} Let $papl$ be the variable capturing the characteristic path length as per constraint~(\ref{papl}).
A lower bound $apl^L$ and an upper bound $apl^U$ can be specified on the characteristic path length via the constraints
\begin{equation}
apl^L \leq papl + sapl^{-} - sapl^{+} \leq apl^U \label{spec-apl},
\end{equation}
where $sapl^{-}$, $sapl^{+} \geq 0$ are a pair of slack variables to be minimized.\qed
\end{example}
\begin{example}
{\bf (Specification of Characteristic Path Length).} Let $pcpl$ be the variable capturing the characteristic path length as per constraint~(\ref{pcpl}).
A lower bound $cpl^L$ and an upper bound $cpl^U$ can be specified on the characteristic path length via the constraints
\begin{equation}
cpl^L \leq pcpl + scpl^{-} - scpl^{+} \leq cpl^U \label{spec-cpl},
\end{equation}
where $scpl^{-}$, $scpl^{+} \geq 0$ are a pair of slack variables to be minimized.\qed
\end{example}
\begin{example}
{\bf (Specification of Diameter).} A lower bound $D^L$ and an upper bound $D^U$ can be specified on the network's diameter via the constraints
\begin{eqnarray}
&1+\left(D^L-1\right)\psi_{ij} \leq w_{ij} + sD^{-} - sD^{+} \leq D^U \;\;\; \forall (i,j) \in E& \label{spec-diam1} \\
&\sum\limits_{(i,j) \in E} \psi_{ij} \geq 1& \label{spec-diam2}
\end{eqnarray}
where $\psi_{ij}$ are auxiliary binary variables defined for each $(i,j) \in E$, and $sD^{-}$, $sD^{+} \geq 0$ are a pair of slack variables to be minimized.
Specifying the diameter to a given value $D$ can be achieved by setting $D^L = D^U = D$.\qed
\end{example}
In the above example, constraints~(\ref{spec-diam1}) impose the strict upper bound of $D^U$ for all shortest paths $w_{ij}$ as well as a lower bound of either $D^L$ or $1$ (a trivial value), depending on whether variable $\psi_{ij}$ assumes the value of $1$ or $0$, respectively. At the same time, constraint~(\ref{spec-diam2}) ensures that ``at least one'' of these lower bounds will be $D^L$.
\begin{example}
{\bf (Specification of Closeness Centrality Sequence).}  Let $piclc_i$ be the variables capturing the inverse of the closeness centrality of each node $i \in V$ as per constraint~(\ref{piclc}).
A sequence $\mathbf{clc} = \left\{\left[clc_1^L,clc_1^U\right], \left[clc_2^L,clc_2^U\right], \cdots, \left[clc_N^L,clc_N^U\right] \right\}$ for the closeness centralities of the nodes in $V$ can be specified via the constraint set\footnote{Constraints~(\ref{spec-clc}) correspond to constraints~(\ref{qim1})--(\ref{qim4}) for the case of the applicable subset being the complete node set $V$.}
\begin{equation}
\hspace{-0.23cm}\left.\begin{array}{c}
\hspace{-0.23cm}\left.\begin{array}{l}
piclc_i + siclc_i^{-} - siclc_i^{+} \geq piclc^L - \left(piclc^L - \frac{1}{clc_m^U} \right) q_{im}\\
piclc_i + siclc_i^{-} - siclc_i^{+} \leq piclc^U - \left(piclc^U - \frac{1}{clc_m^L} \right) q_{im}
\end{array}\right\} \forall (i,m) \in V \times V \\
\hspace{-0.23cm}\left.\begin{array}{l}
\;\sum\limits_{i \in V} q_{im} =  1 \;\;\; \forall m \in V \\
\sum\limits_{m \in V} q_{im} =  1 \;\;\; \forall i \in V,
\end{array}\right.
\end{array}\right.\label{spec-clc}
\end{equation}
where $piclc^L = \min\limits_{j\in V} \left\{\frac{1}{clc_j^U}\right\}$ and $piclc^U = \max\limits_{j\in V} \left\{\frac{1}{clc_j^L}\right\}$ are respectively lower and upper limits for variables $piclc_i$,
$q_{im}$ are $q$-variables to indicate the assignment of the $m^\text{th}$ interval from the sequence $\mathbf{clc}$ to the $i^\text{th}$ node, indicating that the closeness centrality of node $i$, $1/piclc_i$, is to attain a value within the range $\left[clc_m^L,clc_m^U\right]$,
and $siclc_i^{-}, siclc_i^{+} \geq 0$ are two sets of slack variables to be minimized.\qed
\end{example}
\begin{example}
{\bf (Specification of Average Degree of Neighbors).} Let $psdn_i$ be the variables capturing the sum of degrees of the neighbors of each node $i \in V$ as per constraints~(\ref{psdn}), $pnnd_q$ be the variable capturing the number of nodes that attain degree $q$ as per constraints~(\ref{pnnd}), and
$z_{qi}$ be the $z$-variables activated if $pd_i \geq q$ as per constraints~(\ref{pnnd}).
A lower bound $adn_q^L$ and an upper bound $adn_q^U$ can be specified on the average degree of the neighbors of all nodes that attain degree $q$ via the constraints
\begin{eqnarray}
\left.\begin{array}{rcl}
0 \leq &psdn'_{qi}&          \leq \left(d^U \right)^2 \left( z_{iq} - z_{i\left(q+1\right)}\right) \\
0 \leq &psdn_i - psdn'_{qi}& \leq \left(d^U \right)^2 \left( 1 - z_{iq} + z_{i\left(q+1\right)}\right)
\end{array}\right\}  \; \forall i \in V \label{spec-adn1} \\
adn_q^L \; q \; pnnd_q \leq \sum\limits_{i \in V} psdn'_{qi} + sadn_q^{-} - sadn_q^{+} \leq adn_q^U \; q \; pnnd_q \label{spec-adn2},
\end{eqnarray}
where $psdn'_{qi}$ are sets of $p'$-variables corresponding to variables $psdn_i$ (defined for $\tilde{p} = 0$ and activated if node $i$ attains degree $q$),
$d^U$ is the maximum possible degree attainable by any node,
and $sadn_q^{-}, sadn_q^{+} \geq 0$ are two sets of slack variables to be minimized.\qed
\end{example}

In the above example, the variables $psdn'_i$ assume the value $psdn'_i = psdn_i$, if $pd_i = q$ (i.e., if both $z_{iq} = 1$ and $z_{i\left(q+1\right)} = 0$), while otherwise they assume the value $psdn'_i = 0$.
This facilitates the convenient expression of the sum of degrees of all neighbors of nodes that attain a given degree $q$ as the sum of variables $psdn'_i$ over all nodes $i \in V$, irrespective of whether a node $i$ attains degree $q$ or not.
The sum is utilized in constraint~(\ref{spec-adn2}) that imposes the required specification of average degree of the neighbors as a fractional property involving the ratio of this sum over the total number of neighbors of the nodes that attain degree $q$ (which is equal to the product of $q$ times $pnnd_q$).


\subsection{Removing Isomorphism}\label{symmbreak}
A prominent characteristic of the NGP is the potential isomorphism of its solutions;  that is, the admittance of multiple solutions that correspond to the same configuration up to permutation (relabeling) of the node set $V$ and are, thus, practically indistinguishable to each other.
For example, consider a network with $N = 3$ nodes and two configurations $\mathbf{x} = \left(x_{12},x_{13},x_{23}\right)$ such that $\mathbf{x}_A = \left(1,1,0\right)$ and $\mathbf{x}_B = \left(1,0,1\right)$.
It is easy to see that both configurations correspond to the same network, a simple open triplet, the difference being the (arbitrary) label of its central node (node $1$ in $\mathbf{x}_A$ and node $2$ in $\mathbf{x}_B$).
It is not uncommon for networks with as few as ten nodes to have hundreds of isomorphic counterparts.

Isomorphism constitutes a major challenge in addressing the problem via an MILP approach, as the associated branch-and-bound-based search process would be subject to significant burden in exploring unnecessarily enlarged feasible spaces.
Furthermore, the attempt to identify and record multiple distinct solutions will be challenged by the continued recurrence of solutions that are isomorphic to the ones already encountered.
It is therefore desirable to remove all or part of the isomorphism by introducing constraints that will restrict the search to only those solutions that are ``nominal,'' according to some appropriate criterion, while of course ensuring that the full space of distinct, non-isomorphic configurations remains attainable (feasible).
As the goal is essentially to break the symmetry among the nodes in set $V$, the criterion can be based on a node-specific property.
We can then consider that a solution is nominal only if the nodes $i \in V$ attain values for this property, $p_i$, in a sorted (w.l.o.g., non-increasing) fashion:
\begin{equation}
p_i - p_{i+1} \geq 0 \;\;\;\; \forall i \in V \setminus \left\{N\right\}.  \label{symbreak1}
\end{equation}

The degree property is most often the obvious choice of property upon which to reduce isomorphism, since it is easily expressed as a linear sum of the $x$-variables (the formulation's basic variables).
However, constraints~(\ref{symbreak1}) may not suffice to fully eliminate isomorphism due to the lack of a specification that is adequately restrictive (e.g., the specification of a degree sequence with two or more of its elements being equal).
To further break ties among nodes that attain equal values of the ``primary'' property, $p_i$, some ``secondary'' property, $pp_i$, must be considered in tandem.
More specifically, in addition to~(\ref{symbreak1}), we utilize the constraints
\begin{equation}
U \left( p_i - p_{i+1} \right) + pp_i - pp_{i+1} \geq 0 \;\;\;\; \forall i \in V \setminus \left\{N\right\}, \label{symbreak2}
\end{equation}
where $U$ is an appropriate ``big-M'' coefficient that will allow the variables $pp_i$ and $pp_{i+1}$ to vary freely for even a small, yet ``sufficiently'' strictly positive, difference $p_i - p_{i+1}$.
For the case of the the primary symmetry-breaking property $p$ being the degree $p_d$, which is an integral property, the value $U = pp_{i+1}^U - pp_i^L$ would suffice to achieve the desired behavior.

Tertiary and additional, lower-level symmetry-breaking criteria could be imposed in a similar fashion so as to reinforce the attempt to alleviate isomorphism.
However, in most applications there would be little or no residual symmetry to justify the associated increase in model size or the potential scaling challenges due to the introduction of increasingly large coefficients $U^2$, $U^3$, etc.

We remark that breaking symmetry is greatly simplified in the special case where a complete sequence is specified for the primary property (e.g., a full degree sequence).  Let $\mathbf{\Pi} = \left\{\pi_1,\pi_2,\cdots,\pi_N\right\}$ be the desired sequence, which w.l.o.g. we consider to be sorted in non-increasing order.
Primary symmetry-breaking can be achieved by directly specifying that $\pi_i$ be assigned to node $i$, for all nodes $i \in V$.
The secondary symmetry-breaking criterion could then be considered only for those successive node pairs $i$ and $i+1$ that we know a-priori would attain the same degree.
By doing so, not only can symmetry-breaking be imposed without the need for ``big-M'' coefficients, but also the sequence is explicitly assigned without the need for $q$-variables.
Constraints~(\ref{symbreaklite1}) and~(\ref{symbreaklite2}) apply.
The first constraint is a specification of a linear property to an exact value (see constraint~(\ref{lp_bounds}) with $\pi^L = \pi^U \rightarrow \pi_i$) and the second constraint selectively applies the second criterion wherever necessary.
\begin{eqnarray}
&\pi_i \leq p_i \leq \pi_i& \;\;\;\;\; \forall i \in V  \label{symbreaklite1} \\
&pp_i \geq pp_{i+1}&        \;\;\;\;\; \forall i \in V \setminus \left\{N\right\} : \left\{\pi_i = \pi_{i+1}\right\} \label{symbreaklite2}
\end{eqnarray}

\begin{example}
{\bf (Symmetry-breaking based on Degrees and Clustering Coefficients).}  Let $pd_i$ be the variables capturing the degree of each node $i \in V$ as per constraint~(\ref{pd}), and let $pcc_i$ be the variables capturing the (local) clustering coefficient of each node $i \in V$ as per constraints~(\ref{pcc}) or~(\ref{pcc_fix}).
Breaking symmetry among nodes based primarily on their degrees and secondarily on their (local) clustering coefficients can be achieved via the constraint set
\begin{equation}
\left.\begin{array}{r}
pd_i - pd_{i+1} \geq 0 \\
 pd_i - pd_{i+1} + pcc_i - pcc_{i+1} \geq 0
\end{array} \right\}  \forall i \in V \setminus \left\{N\right\},\label{symbreak_deg_cc}
\end{equation}
where the implied coefficient $U=1$ reflects the fact that two clustering coefficients can at most differ by one.
However, when a (non-increasing) degree sequence has been specified as via constraints~(\ref{spec-degseq}), then it suffices to use the simpler constraint set
\begin{equation}
pcc_i \geq pcc_{i+1} \;\;\;\;\; \forall i \in V \setminus \left\{N\right\} : \left\{d_i = d_{i+1}\right\}. \label{symbreak_deg_cc_fix} \qed
\end{equation}
\end{example}

\begin{example}
{\bf (Symmetry-breaking based on Degrees and Shortest-Path Distances).}  Let $pd_i$ be the variables capturing the degree of each node $i \in V$ as per constraint~(\ref{pd}), and let $w_{ij}$ be the variables capturing the shortest-path distance between two nodes $(i,j) \in E$ as per the discussion of Section~\ref{shortestpaths}.
Breaking symmetry among nodes based primarily on their degrees and secondarily on their shortest-path distances from the last node, $i = N$, can be achieved via the constraint set
\begin{equation}
\begin{array}{rl}
pd_i - pd_{i+1} \geq 0 & \;\; \forall i \in V \setminus \left\{N\right\} \\
U \left( pd_i - pd_{i+1} \right) + w_{iN} - w_{\left(i+1\right)N} \geq 0 & \;\; \forall i \in V \setminus \left\{N-1,N\right\},
\end{array} \label{symbreak_deg_wN}
\end{equation}
where $U=N-2$ is the appropriate ``big-M'' coefficient representing the maximum possible difference between two shortest paths in a graph of $N$ nodes.
However, when a (non-increasing) degree sequence has been specified as via constraints~(\ref{spec-degseq}), then it suffices to use the simpler constraint set
\begin{equation}
w_{iN} \geq w_{\left(i+1\right)N} \;\;\;\;\;\; \forall i \in V \setminus \left\{N-1,N\right\} : \left\{d_i = d_{i+1}\right\}. \label{symbreak_deg_wN_fix} \qed
\end{equation}
\end{example}

We remark that, in the above example, any node could have been selected as the reference node for such symmetry breaking purposes.  Yet, in all cases, it would not have been possible to differentiate among the neighbors of the reference node via this set of constraints, since all neighbors feature the same shortest distance of $1$ from the reference node.
Our choice of the last node, $i = N$, as the reference node is thus justified by the fact that it connects to the smallest amount of neighbors (given that the primary breaking of symmetry is done via degrees), which in turn translates to a formulation with as minimal residual isomorphism as possible.

\begin{example}
{\bf (Symmetry-breaking based on Degrees and Sums of Degrees of Neighbors).}  Let $pd_i$ be the variables capturing the degree of each node $i \in V$ as per constraint~(\ref{pd}), and let $psdn_i$ be the variables capturing the sum of degrees of the neighbors of each node $i \in V$ as per constraints~(\ref{psdn}).
Breaking symmetry among nodes based primarily on their degrees and secondarily on the sums of degrees of their neighbors can be achieved via the constraint set
\begin{equation}
\left.\begin{array}{r}
pd_i - pd_{i+1} \geq 0 \\
U \left( pd_i - pd_{i+1} \right) + psdn_i - psdn_{i+1} \geq 0
\end{array} \right\}  \forall i \in V \setminus \left\{N\right\},  \label{symbreak_deg_sdn}\\
\end{equation}
where $U=\left(d^U \right)^2$ is the appropriate ``big-M'' coefficient representing the maximum possible difference between the sums of degrees of neighbors of two given nodes.\footnote{Consider the first node being of maximum degree and connecting to $d^U$ nodes that all attain maximum degree themselves (e.g., as part of a clique), and consider the second node being an isolated node (sum of degrees of neighbors being $0$ in this case).}
However, when a (non-increasing) degree sequence has been specified as via constraints~(\ref{spec-degseq}), then it suffices to use the simpler constraint set
\begin{equation}
psdn_i \geq psdn_{i+1} \;\;\;\;\; \forall i \in V \setminus \left\{N\right\} : \left\{d_i = d_{i+1}\right\}. \label{symbreak_deg_sdn_fix} \qed
\end{equation}
\end{example}

\begin{example}
{\bf (Symmetry-breaking based on Degrees and Closeness Centralities).}  Let $pd_i$ be the variables capturing the degree of each node $i \in V$ as per constraint~(\ref{pd}), and let $piclc_i$ be the variables capturing the (inverse of) closeness centrality of each node $i \in V$ as per constraints~(\ref{piclc}).
Breaking symmetry among nodes based primarily on their degrees and secondarily on their closeness centralities can be achieved via the constraint set
\begin{equation}
\left.\begin{array}{r}
pd_i - pd_{i+1} \geq 0 \\
U \left( pd_i - pd_{i+1} \right) + piclc_{i+1} - piclc_{i} \geq 0
\end{array} \right\}  \forall i \in V \setminus \left\{N\right\},  \label{symbreak_deg_clc}\\
\end{equation}
where $U = \frac{N}{2}-1$ is the appropriate ``big-M'' coefficient representing the maximum possible difference between the inverse closeness centralities of two nodes in any (connected) network,\footnote{The value $\frac{N}{2}-1$ corresponds to the difference between the inverse closeness centrality of a terminal node in a path graph ($cl.c. = \frac{2}{N}$) and that of a node in a clique ($cl.c. = 1$).} and $d^L$, $d^U$ are respectively the minimum and maximum possible degrees attainable by any node.
However, when a (non-increasing) degree sequence has been specified as via constraints~(\ref{spec-degseq}), then it suffices to use the simpler constraint set
\begin{equation}
piclc_i \leq piclc_{i+1} \;\;\;\;\; \forall i \in V \setminus \left\{N\right\} : \left\{d_i = d_{i+1}\right\}. \label{symbreak_deg_clc_fix} \qed
\end{equation}
\end{example}

\section{Illustrative Applications}\label{studies}
In this section we discuss the development of complete formulations for a number of informative applications.
We demonstrate the modularity of the approach by utilizing constraint sets introduced previously in Sections~\ref{popular} and~\ref{enhancements}.
All case studies were solved on a Win7 machine with an Intel i5-3230M (2.60GHz) processor and 3GB RAM. Gurobi~5.6~\cite{GUROBI} was used as the MILP optimization solver, and all runs were explicitly limited to utilize a single CPU core.
Aside from customized branching prioritization for binary variables according to the earlier discussion, default solver options were used throughout this study.
The node labels in all figures represent the degrees of the depicted nodes.

\subsection{Tuning Network Connectivity}
Two primary characteristics of a network in terms of the connectivity among its nodes are the degree distribution and the level of clustering.
As elaborated in the introduction, there exist many well-known algorithms for constructing networks with a specified degree sequence.
However, none of these algorithms can guarantee that the constructed networks adhere also to pre-specified clustering statistics, which can in general vary significantly among networks of the given degree sequence.

In fact, networks generated by approaches that randomly connect vertices tend to have a low number of triangles and, thus, low clustering coefficients.
In contrast, many networks that develop in the real world feature larger clustering coefficients.
Take social networks, for instance, where two randomly chosen friends (neighbors) of a person (node) are also more likely to know each other than two randomly chosen people in the network.
Furthermore, networks whose nodes have an underlying spatial meaning (e.g., power-grid networks) tend to develop more edges between nodes that are closer to each other, leading to the presence of more triangles than what would be expected by mere chance.
To that end, it is sometimes necessary to be able to control the statistics of both degrees and clustering of the nodes in the constructed network.

In this application, we consider specifying the degree sequence in conjunction with the network's average and global clustering coefficients.
Let $\mathbf{d} = \left\{d_1, d_2, \cdots, d_N\right\}$ be a specification for the degree sequence\footnote{For the clustering coefficients to be well-defined, it should hold that $d_N \geq 2$.} (w.l.o.g., we assume $d_i \geq d_{i+1} \; \forall i \in V \setminus \left\{N\right\}$), and let $\left[acc^L,acc^U\right]$ and $\left[gcc^L,gcc^U\right]$ be respectively the acceptable bounds for the average and global clustering coefficient.
After picking the appropriate sets of constraints from examples presented previously in the paper, we formulate the following optimization model:
\begin{eqnarray}
&\min\limits_{\mathbf{x}} & \sum\limits_{i \in V} \left(sd_i^{-} + sd_i^{+}\right) + sacc^{-} + sacc^{+} + sgcc^{-} + sgcc^{+} \label{formclu-obj} \\
& \text{s.t.} & \text{Constraints~(\ref{ytr})} \;\;\;\;\;\;\; \forall (i,j,k) \in T  \nonumber \\
&             & \text{Constraints~(\ref{pd})} \;\;\;\;\;\;\; \forall i \in V  \nonumber \\
&             & \text{Constraints~(\ref{pntr})} \;\;\;\;\;\;\; \forall i \in V  \nonumber \\
&             & \text{Constraints~(\ref{pcc_fix})} \;\;\;\;\;\;\; \forall i \in V  \nonumber \\
&             & \text{Constraint~\,\,(\ref{pacc})} \nonumber \\
&             & \text{Constraint~\,\,(\ref{pgcc_fix})} \nonumber \\
&             & \text{Constraint~\,\,(\ref{spec-degseq})} \nonumber \\
&             & \text{Constraint~\,\,(\ref{spec-acc})} \nonumber \\
&             & \text{Constraint~\,\,(\ref{spec-gcc_fix})} \nonumber \\
&             & \text{Constraints~(\ref{symbreak_deg_cc_fix})} \nonumber \\
&& x_{ij} \in \left\{0,1\right\} \;\;\;\;\;\;\;\;\;\;\;\;\;\;\;\; \forall (i,j) \in E \label{formclu-xbinary}\\
&& sd_i^{-},sd_i^{+} \geq 0 \;\;\;\;\;\;\;\;\;\;\;\;\;\; \forall i \in V \label{formclu-slackpos1}\\
&& sacc^{-},sacc^{+},sgcc^{-},sgcc^{+} \geq 0 \label{formclu-slackpos2}.
\end{eqnarray}

Constraints~(\ref{ytr}) through~(\ref{pgcc_fix}) restrict all applicable $p$-variables (variables $pd_i$, $pacc$ and $pgcc$) to their appropriate values (the degree of each node and the two clustering coefficients, respectively), while constraints~(\ref{spec-degseq}), (\ref{spec-acc}) and~(\ref{spec-gcc_fix}) impose the specifications.
Obviously, if no specification is desired for either of the two clustering coefficients, the corresponding constraint~(\ref{spec-acc}) or~(\ref{spec-gcc_fix}) should be omitted (or simply ``disabled'' by setting $\left[0,1\right]$ as the allowable bounds).
Constraints~(\ref{spec-degseq}) are also used in conjunction with the ordered degree sequence specification to reduce isomorphism.
Symmetry is primarily broken based on node's degrees, while residual symmetry is broken based on local clustering coefficients.
Finally, constraints~(\ref{formclu-xbinary}) declare the binarity of the $x$-variables, while constraints~(\ref{formclu-slackpos1}) and~(\ref{formclu-slackpos2}) declare the non-negativity of the $s$-variables.
The objective~(\ref{formclu-obj}) minimizes the sum of every slack variable associated with a specification.
We remark that although a number of $y$-, $p$- and $s$-variables participate in the formulation, we only list $\mathbf{x}$ under the $\min$ operator to highlight the fact that the $x$-variables suffice to fully define the solution (i.e., the values of all other variables can be uniquely defined once values for $\mathbf{x}$ are determined).

\begin{figure}[!ht]
\captionsetup[subfigure]{labelformat=empty}
\centering
\subfloat[]{
\includegraphics[totalheight=0.27\textheight,trim=0 2.0cm 0 2.5cm,clip=true]{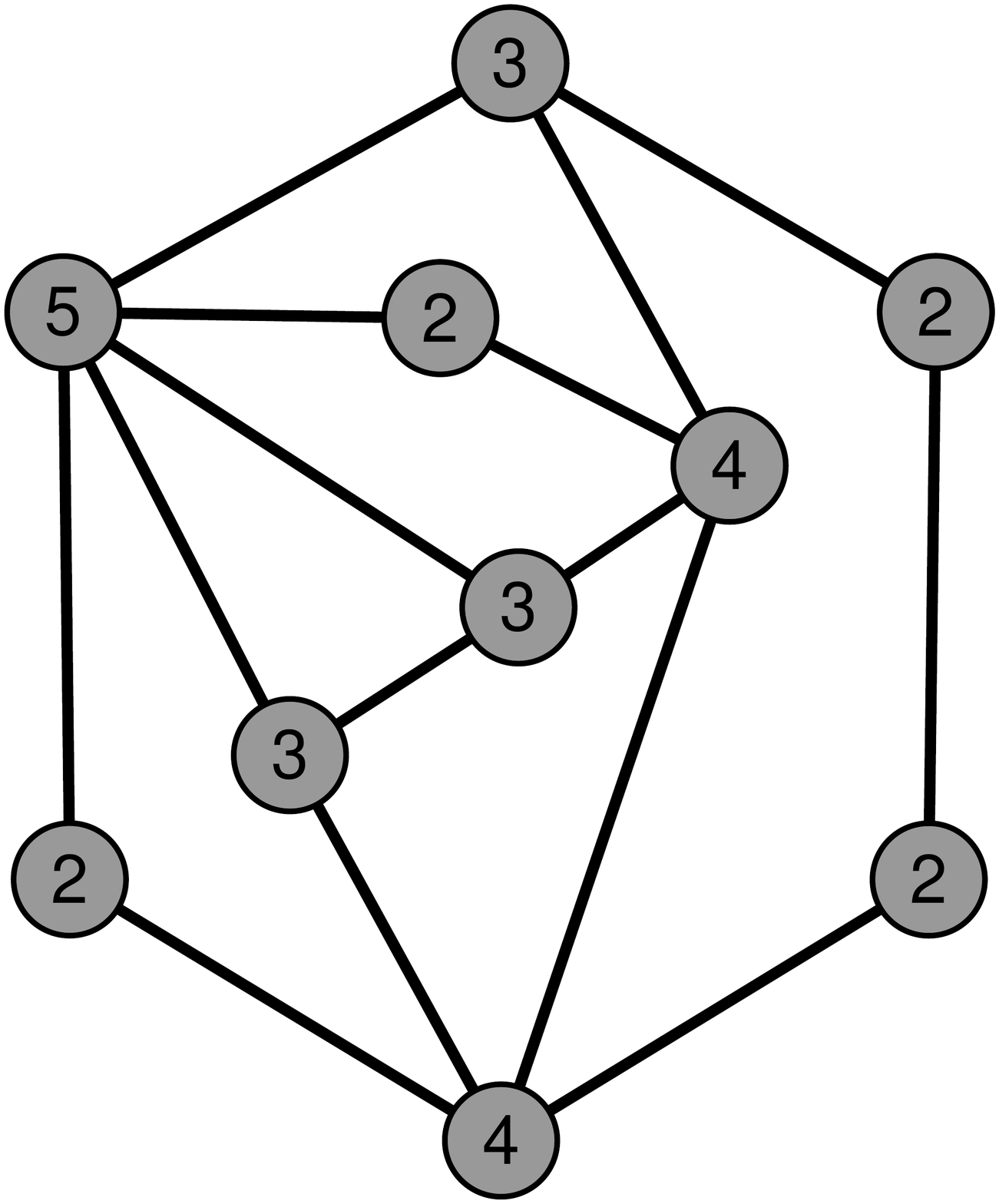}
}
\subfloat[]{
\includegraphics[totalheight=0.30\textheight,trim=0 5.0cm 0 3.5cm,clip=true]{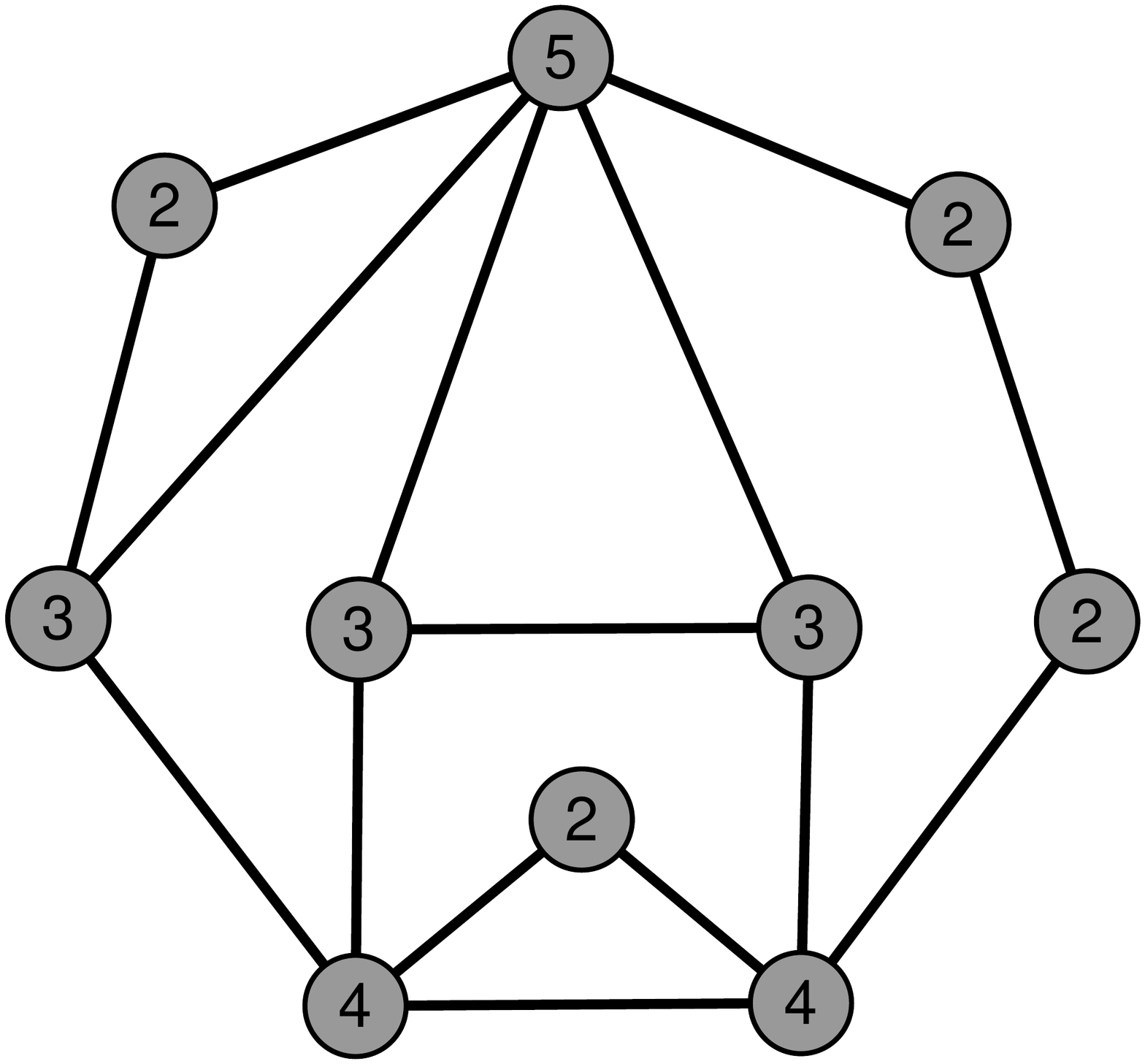}
}
\subfloat[]{
\includegraphics[totalheight=0.26\textheight,trim=0 1.5cm 0 1.5cm,clip=true]{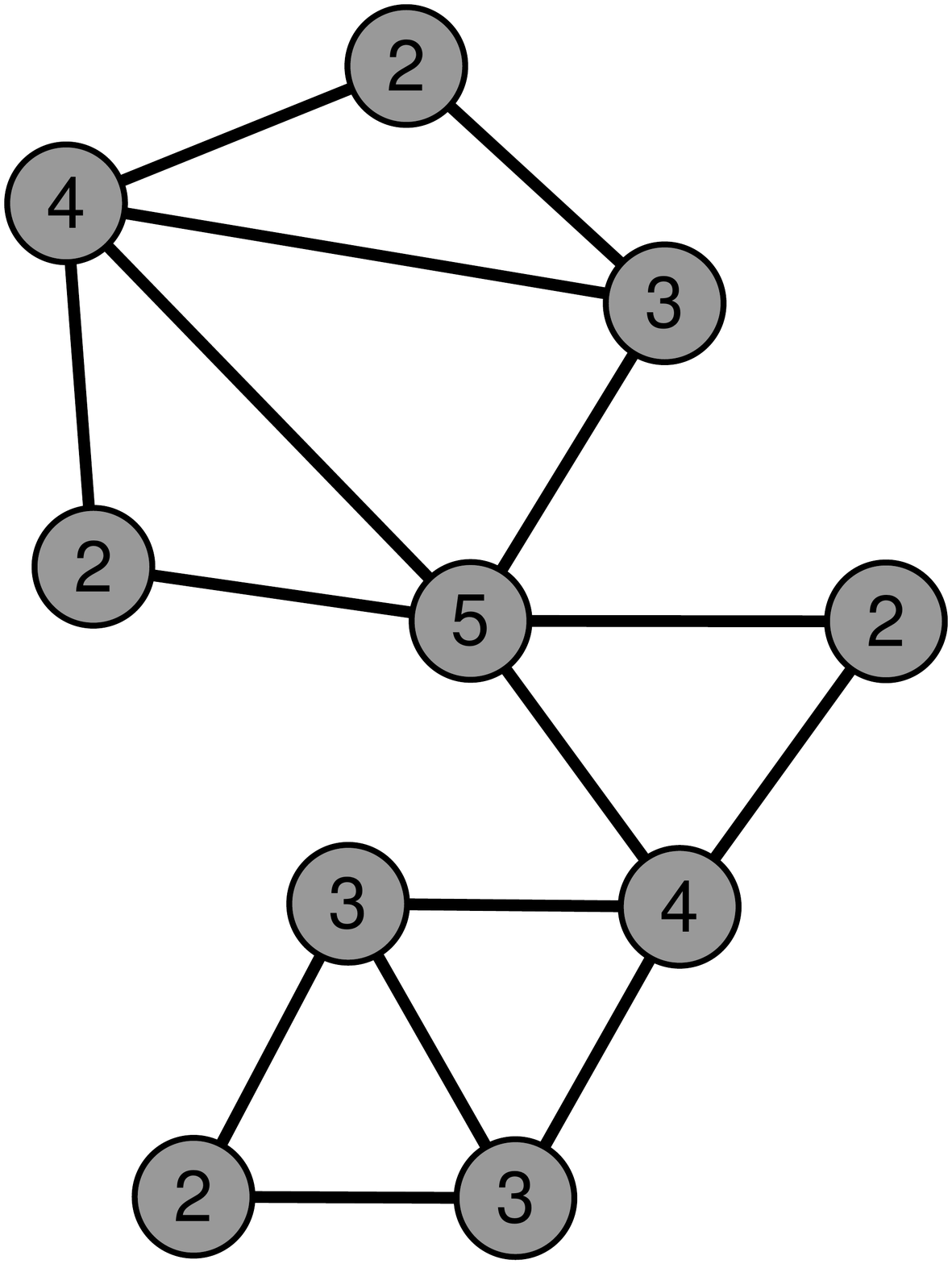}
}
\caption{Representative optimal solutions to the three specifications from Case Study~\ref{case1}.  From left to right, the networks possess the same degree sequence but are respectively tuned to low ($acc = 0.077, gcc = 0.086$), medium ($acc = 0.353, gcc = 0.257$) and high ($acc = 0.713, gcc = 0.514$) level of clustering.}
\label{fig:case1}
\end{figure}

Case Study~\ref{case1} demonstrates the use of this formulation.
We provide a specific (randomly chosen) degree sequence, and we desire to construct networks with various values for the average and the global clustering coefficient.
To that end, we provide three ranges (``low,'' ``medium'' and ``high'') from which we would require both coefficients to attain values.

\begin{casestudy}\label{case1}
Let $N=10$ nodes, and let $\mathbf{d} =\left\{5, 4, 4, 3, 3, 3, 2, 2, 2, 2\right\}$ be a specification for the degree sequence.
Let also three separate specifications with $\left[0.000,0.250\right]$, $\left[0.250,0.500\right]$
and $\left[0.500,0.750\right]$ as the allowable bounds for the average and global
clustering coefficients (bounds common to both metrics).  %
Figure~\ref{fig:case1} depicts representative solutions.
Each of the three networks was constructed in less than 0.1 $sec$ CPU. \qed
\end{casestudy}

\begin{figure}[!ht]
\captionsetup[subfigure]{labelformat=empty}
\centering
\subfloat[]{
\includegraphics[totalheight=0.27\textheight,trim=0 4.5cm 0 4.5cm,clip=true]{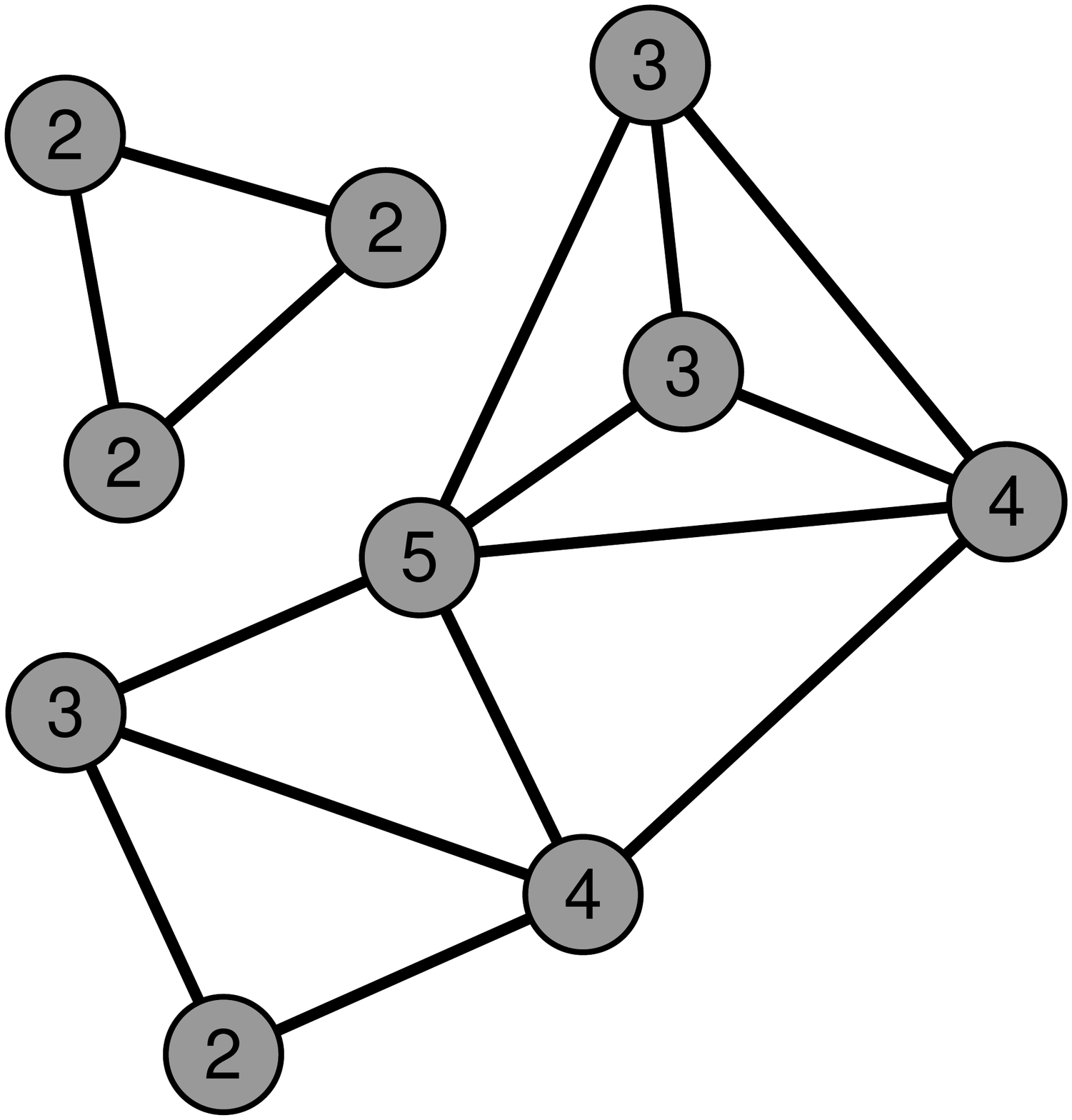}
}
\caption{A common optimal solution to the problems of maximizing the average and global clustering coefficients given a specified degree sequence. The network possesses maximum values for the average and global clustering coefficients, $acc = 0.833$ and $gcc = 0.686$, respectively.}
\label{fig:case2}
\end{figure}

Instead of trying to locate network realizations that adhere to a particular specification, a problem cast as a slack-minimization problem, Case Study~\ref{case2} discusses the situation where an application-specific objective function is available.
In particular, we seek to identify the network with maximum clustering coefficient (either average or global) that adheres to the specified degree sequence.
The applicable formulations can be derived from the previous one per the following steps: (i) remove all $s$-variables (i.e., set them to a fixed value of zero), (ii) disable the specifications for the two coefficients (i.e., set $acc^L=gcc^L=0$ and $acc^U=gcc^U=1$ in constraints~(\ref{spec-acc}) and~(\ref{spec-gcc_fix}), or completely remove these constraints from the formulation), and (iii) replace the objective with either $\left\{\max\limits_{\mathbf{x}} pacc\right\}$ or $\left\{\max\limits_{\mathbf{x}} pgcc\right\}$, depending on which of the two coefficients we seek to maximize.

\begin{casestudy}\label{case2}
Let $N=10$ nodes, and let $\mathbf{d} =\left\{5, 4, 4, 3, 3, 3, 2, 2, 2, 2\right\}$ be a specification for the degree sequence.
We seek the networks that maximize the average clustering coefficient as well as the networks that maximize the global clustering coefficient.
Figure~\ref{fig:case2} depicts a common solution that happens to maximize both metrics at the same time.
The network was obtained in about 1 $sec$ CPU. \qed
\end{casestudy}

We observe that Case Study~\ref{case2} resulted in a network that features two disconnected components.  The optimizer was able to discern that, in this way, ``more clustered'' communities can be formed while at the same time adhering to the overall specified degree sequence.

\subsection{Tuning Network Spread}
Beyond local connectivity in the network, which can be controlled through properties such as the degrees and clustering coefficients, another informative network quality to consider during
construction is the network's {\em long-range spread}.
The spread is especially important in networks where the edges represent conduits of information flow, such as the case of communication networks.
It relates to the relative size of shortest paths across the network and can thus be characterized by a number of relevant properties, most notably the characteristic path length (cpl) and the network diameter (i.e., the median and maximum shortest-paths, respectively).

It is important to note that the network spread may sometimes correlate with node connectivity, but tuning connectivity does not suffice for also tuning the spread.
Randomly constructed networks typically have a small number of triangles (low clustering) and low spread.
On the other hand, small-world networks also exhibit low spread, but are known to feature a large number of triangles.
In general, networks with different possible combinations of local connectivity and spread exist, and it is important to be able to independently control both qualities during network construction.

To that end, in this application we study the construction of networks that adhere to a specific degree sequence and at the same time possess pre-specified values for the two network spread metrics, cpl and diameter.
Let $\mathbf{d} = \left\{d_1, d_2, \cdots, d_N\right\}$ be a specification for the degree sequence (w.l.o.g., we assume $d_i \geq d_{i+1} \; \forall i \in V \setminus \left\{N\right\}$), and let $\left[cpl^L,cpl^U\right]$ and $\left[D^L,D^U\right]$ be respectively the acceptable bounds for the network's cpl and diameter.
The following formulation applies:
\begin{eqnarray}
&\min\limits_{\mathbf{x}} & \sum\limits_{i \in V} \left(sd_i^{-} + sd_i^{+}\right) + scpl^{-} + scpl^{+} + sD^{-} + sD^{+} \label{formspread-obj} \\
& \text{s.t.} & \text{Constraints~(\ref{sp-p1})--(\ref{sp-p4}),~(\ref{sp-d1})--(\ref{sp-d4}),~(\ref{sp-o1}),~(\ref{sp-o2})} \;\;\; \forall (i,j) \in E  \nonumber \\
&             & \text{Constraints~(\ref{pd})} \;\;\;\;\;\;\;\;\;\;\;\;\;\, \forall i \in V  \nonumber \\
&             & \text{Constraints~(\ref{pcpl})} \nonumber \\
&             & \text{Constraints~(\ref{spec-degseq})} \nonumber \\
&             & \text{Constraints~(\ref{spec-cpl})} \nonumber \\
&             & \text{Constraints~(\ref{spec-diam1}),~(\ref{spec-diam2})} \nonumber \\
&             & \text{Constraints~(\ref{symbreak_deg_wN_fix})} \nonumber \\
&& x_{ij}, \psi_{ij} \in \left\{0,1\right\} \;\;\;\;\;\;\;\;\;\;\;\;\;\;\; \forall (i,j) \in E \label{formspread-xbinary} \\
&& sd_i^{-},sd_i^{+} \geq 0 \;\;\;\;\;\;\;\;\;\;\;\;\;\;\;\;\;\;\;\; \forall i \in V \label{formspread-slackpos1} \\
&& scpl^{-},scpl^{+},sD^{-},sD^{+} \geq 0 \label{formspread-slackpos2}.
\end{eqnarray}

Constraints~(\ref{sp-p1}) through~(\ref{sp-o2}) correspond to the encoding in $w$-variables of shortest path lengths via the primal/dual formulation analysis presented in Section~\ref{shortestpaths}.
Constraints~(\ref{pd}) and~(\ref{pcpl}) define the applicable $p$-variables for the node degrees and the cpl, while constraints~(\ref{spec-degseq}) and~(\ref{spec-cpl}) impose the corresponding specifications. The network diameter is specified via constraints~(\ref{spec-diam1}) and~(\ref{spec-diam2}).
Obviously, if no specification is desired for either the cpl or the diameter, the corresponding constraints should be omitted (or simply ``disabled'' by setting sufficiently ``loose'' bounds).
Isomorphism is primarily eliminated based on node degrees, while residual symmetry in the problem is broken via constraints~(\ref{symbreak_deg_wN_fix}), which order nodes of the same degree according to their shortest path to the last node (node $i = N$).
Finally, constraints~(\ref{formspread-xbinary}) declare the binarity of the $x$- and $\psi$-variables, while constraints~(\ref{formspread-slackpos1}) and~(\ref{formspread-slackpos2}) declare the non-negativity of the $s$-variables.
The objective~(\ref{formspread-obj}) minimizes the sum of every slack variable associated with a specification.

\begin{figure}[!ht]
\captionsetup[subfigure]{labelformat=empty}
\centering
\subfloat[]{
\includegraphics[totalheight=0.27\textheight,trim=0 4.5cm 0 4.5cm,clip=true]{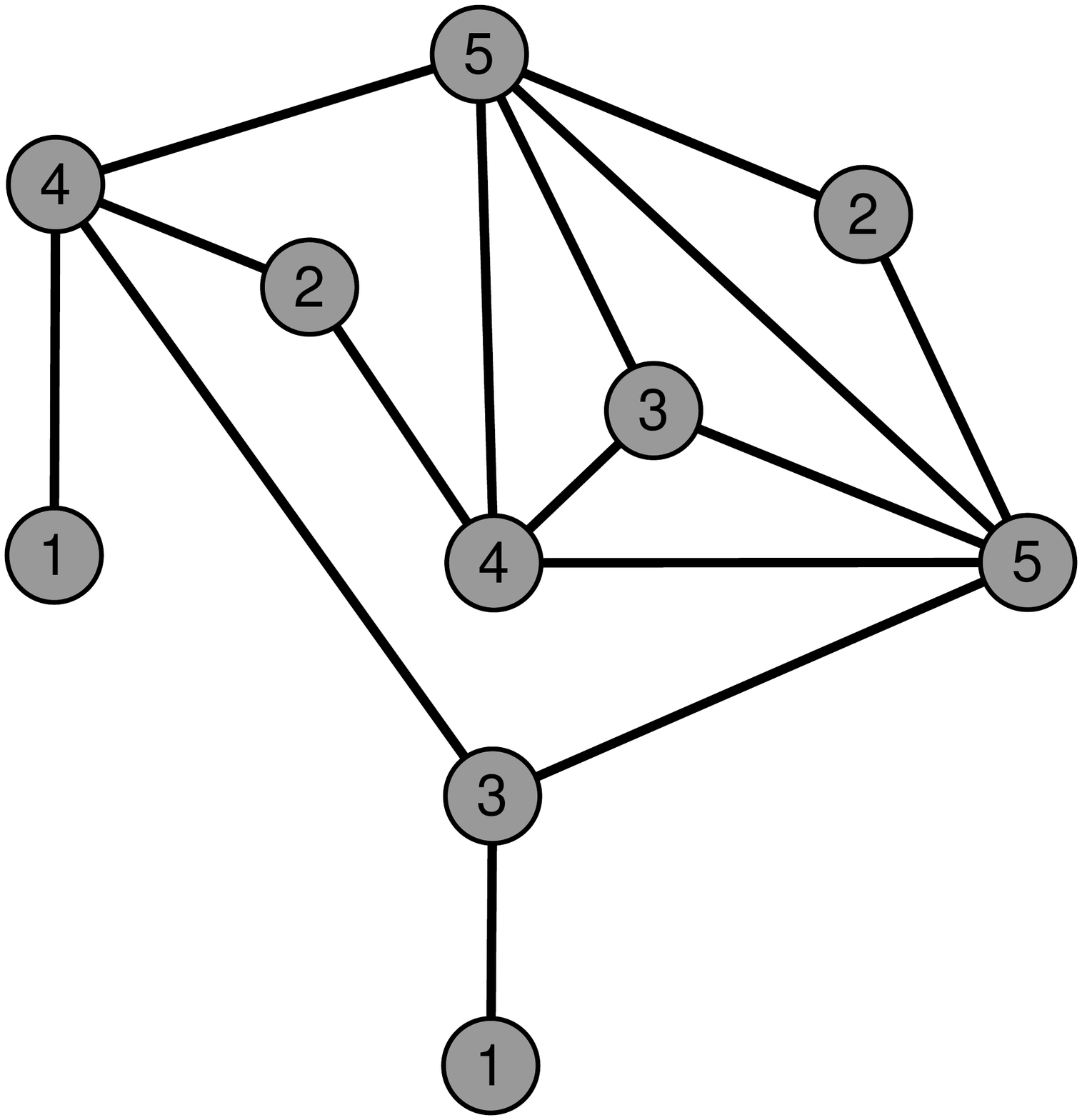}
}
\subfloat[]{
\includegraphics[totalheight=0.27\textheight,trim=0 0.0cm 0 0.0cm,clip=true]{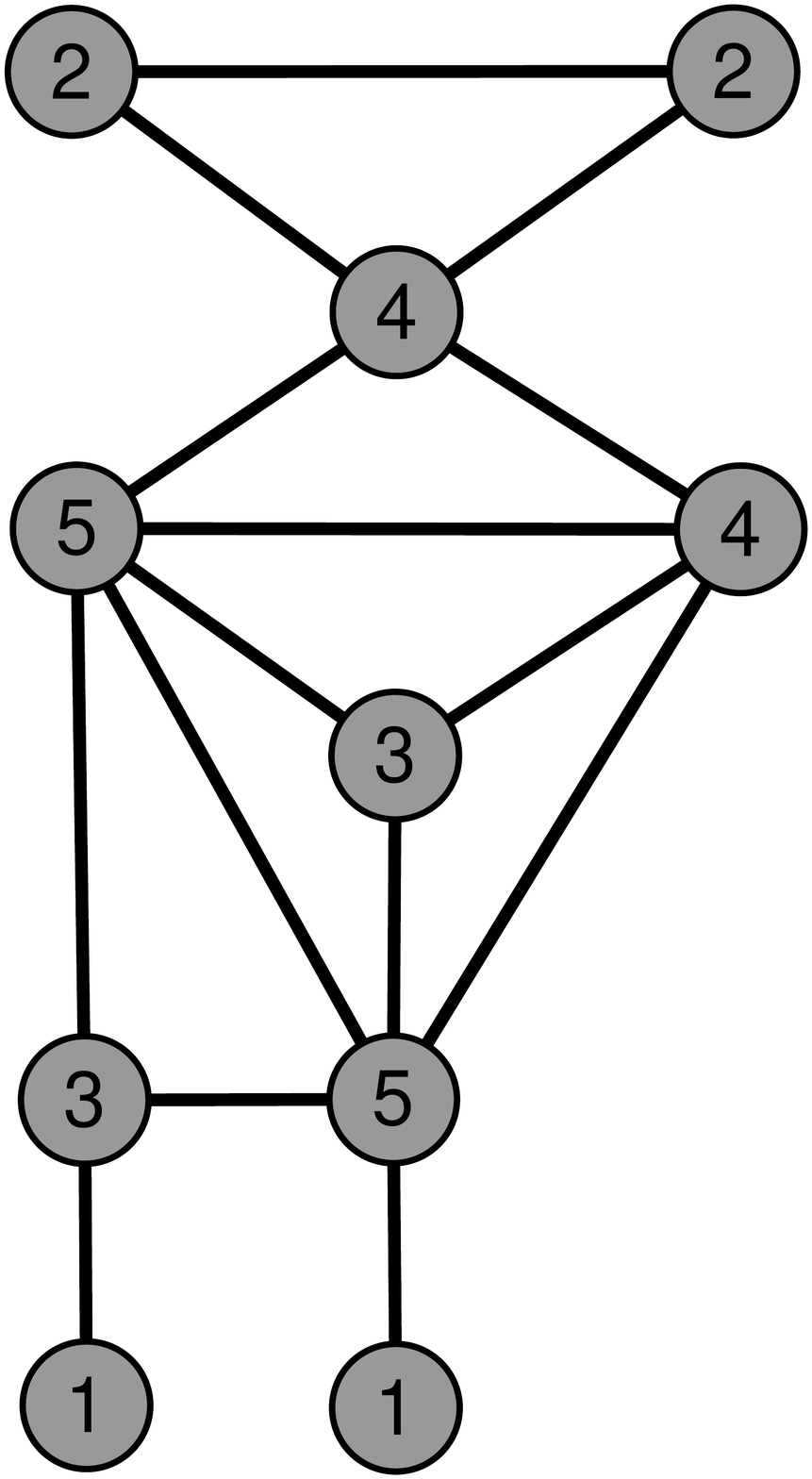}
}
\subfloat[]{
\includegraphics[totalheight=0.27\textheight,trim=0 3.5cm 0 3.5cm,clip=true]{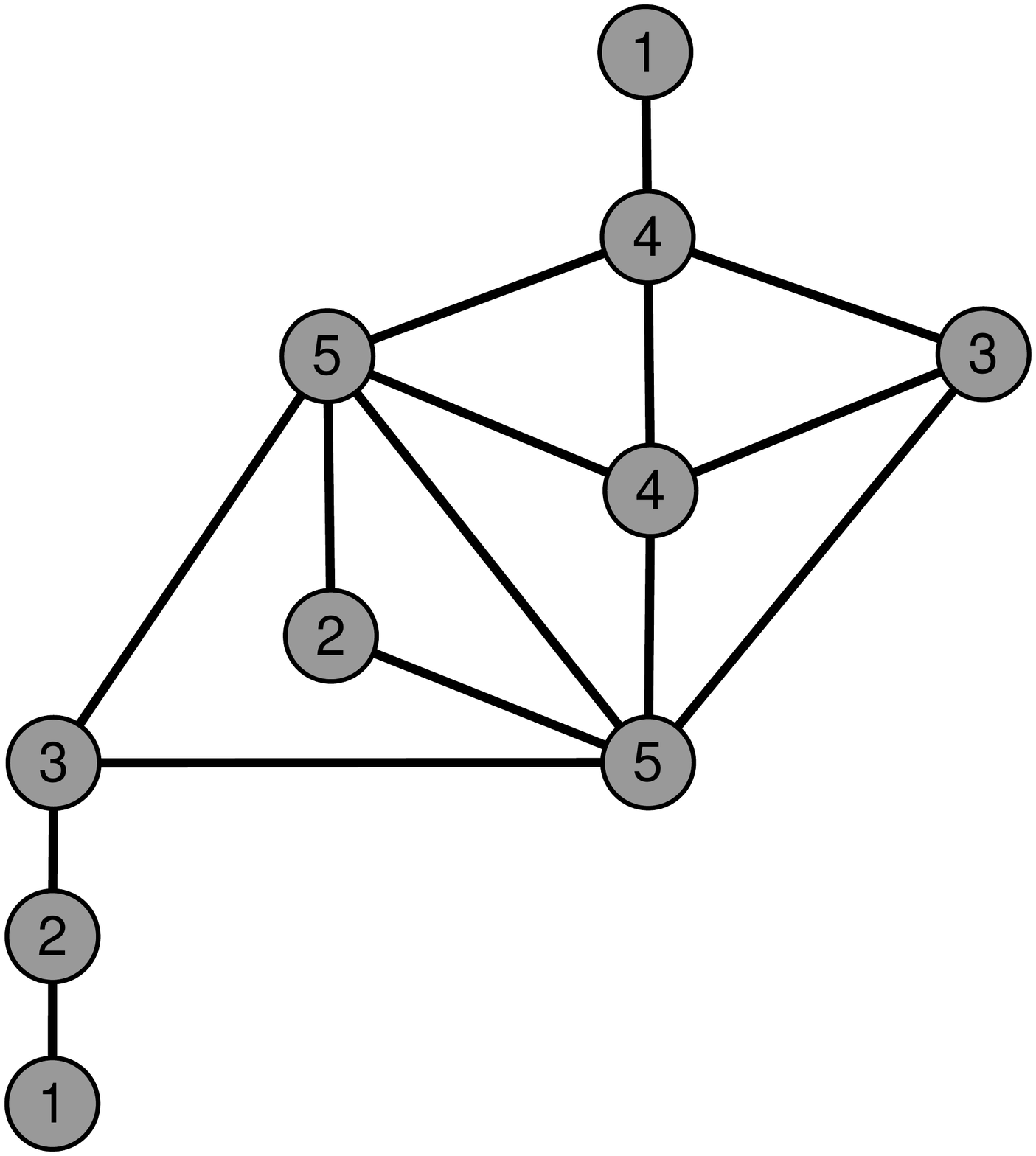}
}
\caption{Representative optimal solutions to the three specifications from Case Study~\ref{case3}.  From left to right, the networks possess diameter values of 3, 4 and 5.  All three networks possess a characteristic path length value of 2.}
\label{fig:case3}
\end{figure}

Case Study~\ref{case3} demonstrates the use of this formulation. We provide a degree sequence as in the previous case studies, and we desire to construct networks that attain various combinations of values for the cpl and the diameter.

\begin{casestudy}\label{case3}
Let $N=10$ nodes, and let $\mathbf{d} =\left\{5, 5, 4, 4, 3, 3, 2, 2, 1, 1\right\}$ be a specification for the degree sequence.
In addition, let a specification of $cpl^L = cpl^U = 2$ for the characteristic path length, and let three separate specifications for the diameter requiring that it attains the values of $D=3$, $D=4$ and $D=5$.
Figure~\ref{fig:case3} depicts representative solutions for the three cases, which were respectively obtained in approximately 20, 23 and 198 $sec$ CPU. \qed
\end{casestudy}

\subsection{Tuning Network Assortativity}
In the previous applications we specified node degrees, clustering information and statistics on the lengths of shortest paths.
In this section we will consider mixing patterns between nodes that form constituent links; that is, correlations between node properties of neighbors.
We specifically focus on capturing the correlations between degrees of neighbors in a network.
A network is said to have assortative mixing if high (low) degree nodes preferentially form links with high (low) degree nodes.
In contrast, disassortative networks are those in which high degree nodes preferentially link with low degree ones, and vice versa.
Empirical studies of real world networks have revealed that social networks tend to be assortative, while technological and biological networks are usually disassortative.
In either case, having the ability to accurately capture degree-degree correlations of real world networks is crucial to understanding the dynamic processes that govern the formation of these networks.

Degree-degree correlations can be quantified in different ways.
We pose the problem such that we specify bounds for the average degree among all neighbors of nodes of a given degree.
For convenience, we further specify an overall minimum and maximum degree for the whole network (which may just be the trivial values of 0 and $N-1$, respectively).
Note that we do not specify the exact degree sequence as in the previous applications.

More formally, let $Q = \{d^L,d^L+1,\cdots,d^U-1,d^U\}$ be the set of possible degrees, where $0 \leq d^L \leq d^U \leq N-1$, and let also $Q^+ = Q \setminus \left\{0\right\}$ defined for notational convenience.
We focus on each possible degree $q \in Q^+$ and require that the average degree of the
neighbors of all nodes that attain degree $q$, denoted as $adn_q$, be
restricted in the range $\left[adn_q^L,adn_q^U\right]$.
Note that we do not provide a specification for $q=0$, as isolated nodes do not
possess any neighbors.
We formulate the following optimization model:
\begin{eqnarray}
&\min\limits_{\mathbf{x}} & \sum\limits_{i \in V} \left(sd_i^{-} + sd_i^{+}\right) + \sum\limits_{q \in Q^+} \left(sadn_q^{-} + sadn_q^{+}\right) \label{formassort-obj} \\
&             & \text{Constraints~(\ref{pd})} \;\;\;\;\;\;\;\;\;\;\;\;\; \forall i \in V  \nonumber \\
&             & \text{Constraints~(\ref{psdn})} \nonumber \\
&             & \text{Constraints~(\ref{pnnd})} \nonumber \\
&             & \text{Constraints~(\ref{spec-deg})} \;\;\;\;\;\;\;\;\;\;\;\;\; \forall i \in V  \nonumber \\
&             & \text{Constraints~(\ref{spec-adn1}),~(\ref{spec-adn2})} \;\;\;\; \forall q \in Q^+ \nonumber \\
&             & \text{Constraints~(\ref{symbreak_deg_sdn})} \nonumber \\
&& pd_1 \geq 1 \label{formassort-nonull}\\
&& x_{ij} \in \left\{0,1\right\} \;\;\;\;\;\;\;\;\;\;\;\;\;\;\;\; \forall (i,j) \in E \label{formassort-xbinary} \\
&& sd_i^{-},sd_i^{+} \geq 0 \;\;\;\;\;\;\;\;\;\;\;\;\;\;\, \forall i \in V \label{formassort-slackpos1} \\
&& sadn_q^{-},sadn_q^{+} \geq 0 \;\;\;\;\;\;\,\forall q \in Q^+ \label{formassort-slackpos2}.
\end{eqnarray}

Constraints~(\ref{pd}),~(\ref{psdn}) and~(\ref{pnnd}) define variables $pd_i$, $psdn_i$ and $pnnd_q$ for the degrees, the sums of neighbors' degrees and the number of nodes that attain a given degree, respectively.
Constraints~(\ref{spec-deg}) impose the degree specification, as it applies to each node independently, and they are only relevant when bounds other than $d^L = 0$ and $d^U = N-1$ are chosen.
Constraints~(\ref{spec-adn1}) and~(\ref{spec-adn2}) impose the specification for the average degree of neighbors of the nodes of each degree $q \in Q^+$.
Constraints~(\ref{symbreak_deg_sdn}) break symmetry, primarily based on degrees, with residual symmetry among nodes that attain the same degree being broken based on the sum of degrees of their neighbors.
Constraint~(\ref{formassort-nonull}) is added to avoid the trivial solution of a null graph, which may otherwise be feasible under certain specifications.
Finally, constraints~(\ref{formassort-xbinary}) through~(\ref{formassort-slackpos2}) declare binarity and positivity of variables.
As before, the objective~(\ref{formassort-obj}) minimizes the sum of the slack variables used in the formulation.
For the offset parameter in constraints~(\ref{pnnd}) we chose the value $\epsilon = 0.01$, which is small enough to not adversely affect the tightness of the LP-relaxations but at the same time large-enough to sufficiently exceed the optimization solver's (default) feasibility tolerance setting of $10^{-6}$.

\begin{figure}[!ht]
\captionsetup[subfigure]{labelformat=empty}
\centering
\subfloat[]{
\includegraphics[totalheight=0.27\textheight,trim=0 5.0cm 0 4.0cm,clip=true]{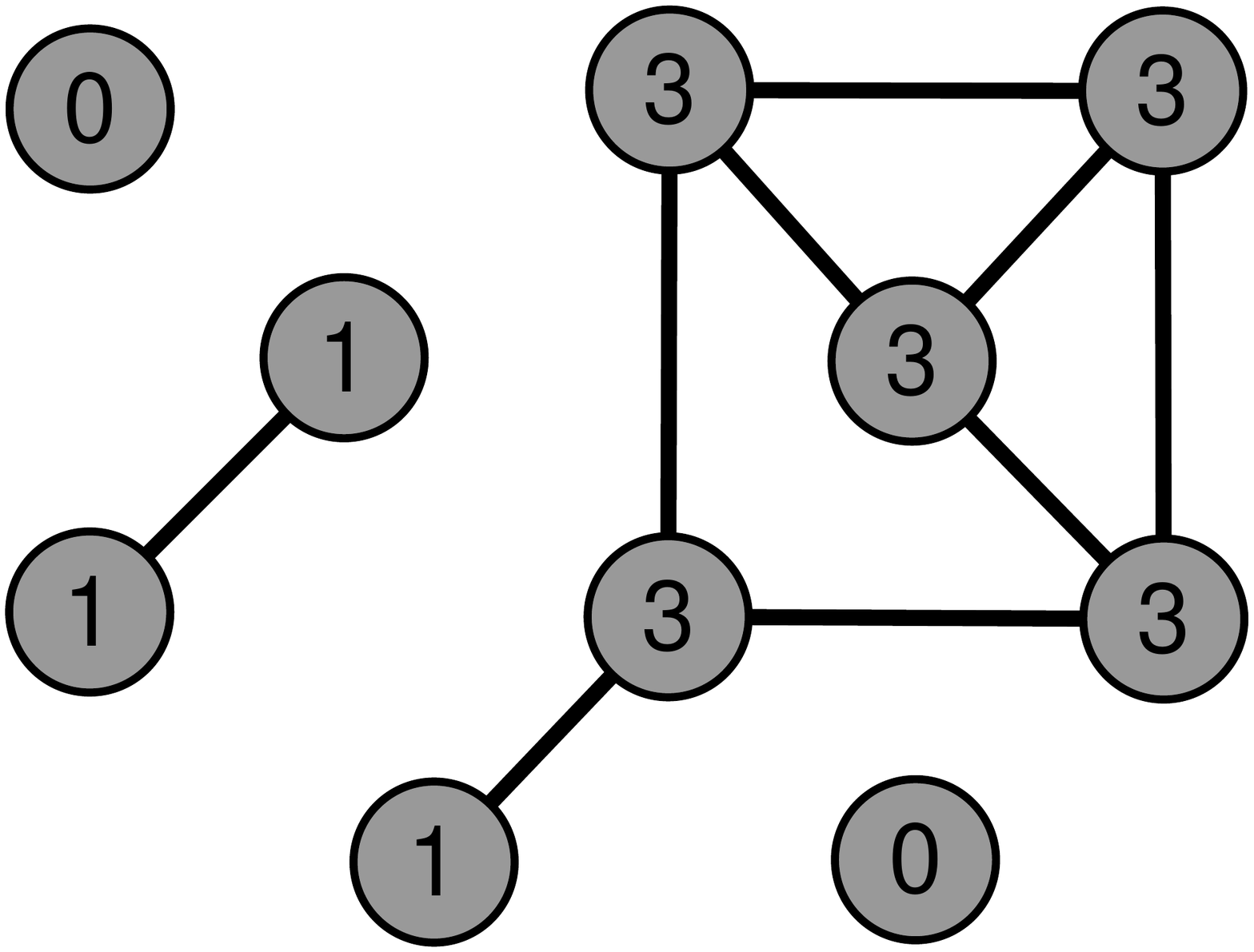}
}
\hspace{1cm}
\subfloat[]{
\includegraphics[totalheight=0.27\textheight,trim=0 5.0cm 0 5.0cm,clip=true]{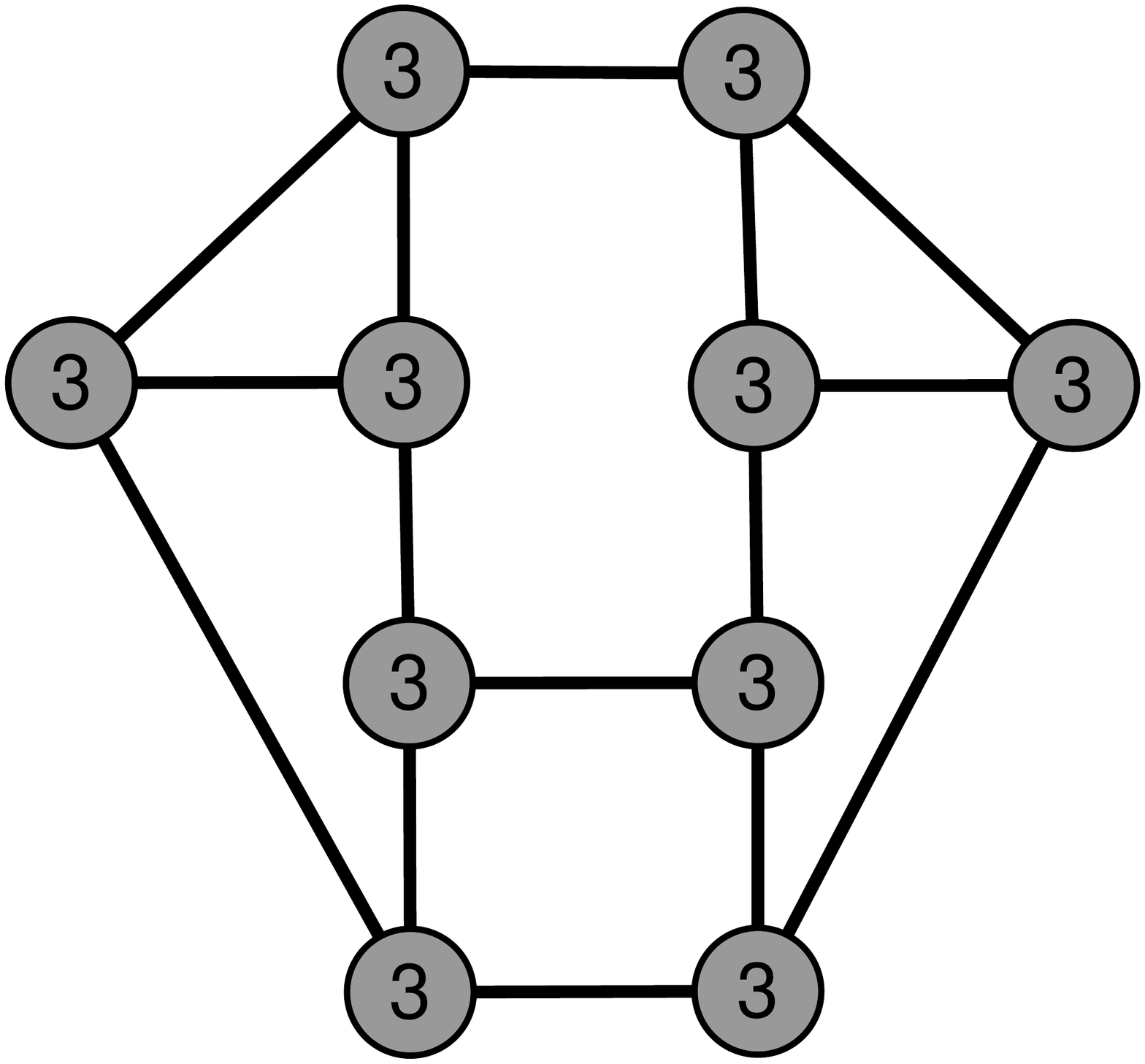}
}\hfil
\subfloat[]{
\includegraphics[totalheight=0.27\textheight,trim=0 5.0cm 0 4.0cm,clip=true]{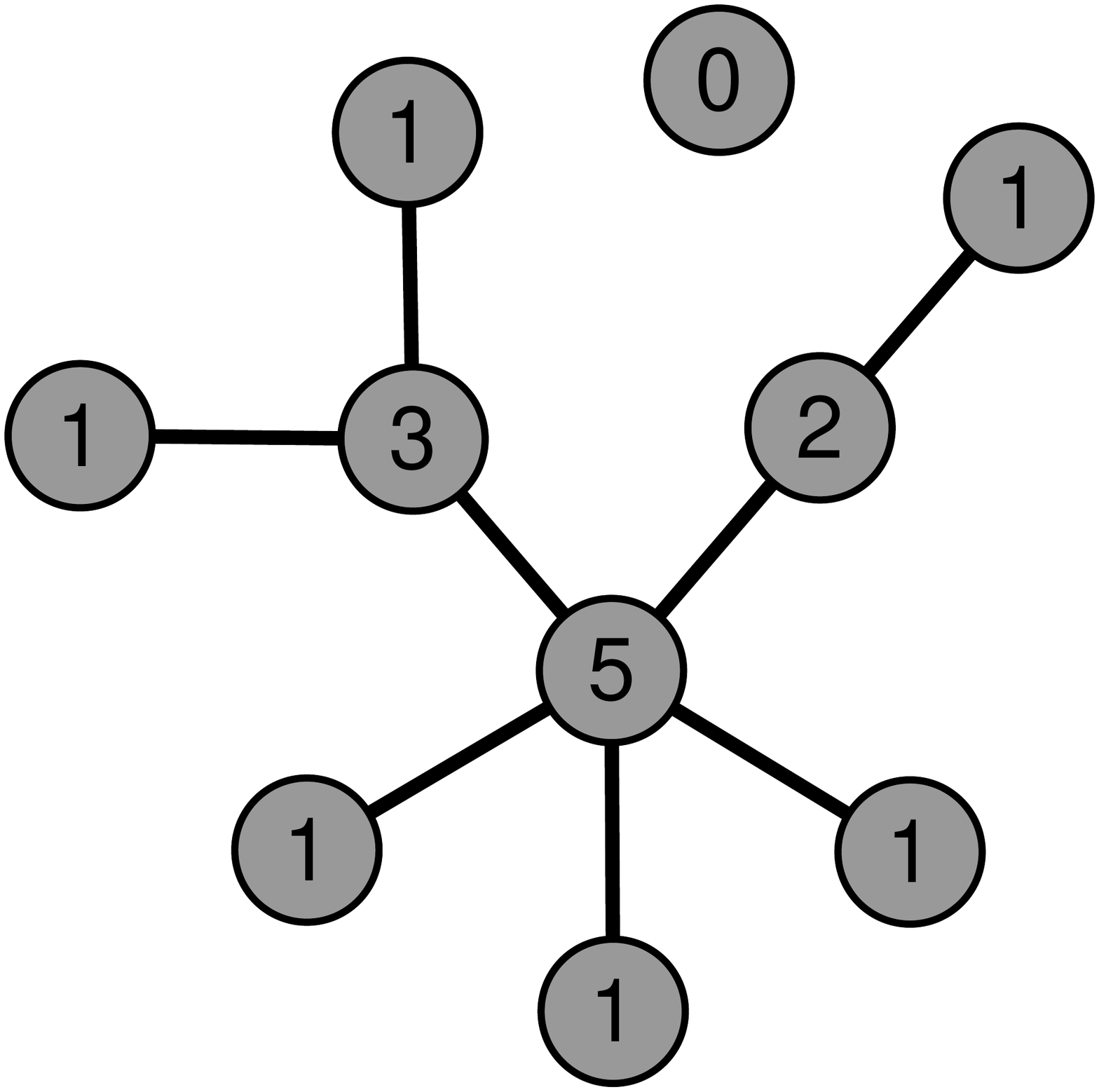}
}
\hspace{1cm}
\subfloat[]{
\includegraphics[totalheight=0.27\textheight,trim=0 5.0cm 0 5.0cm,clip=true]{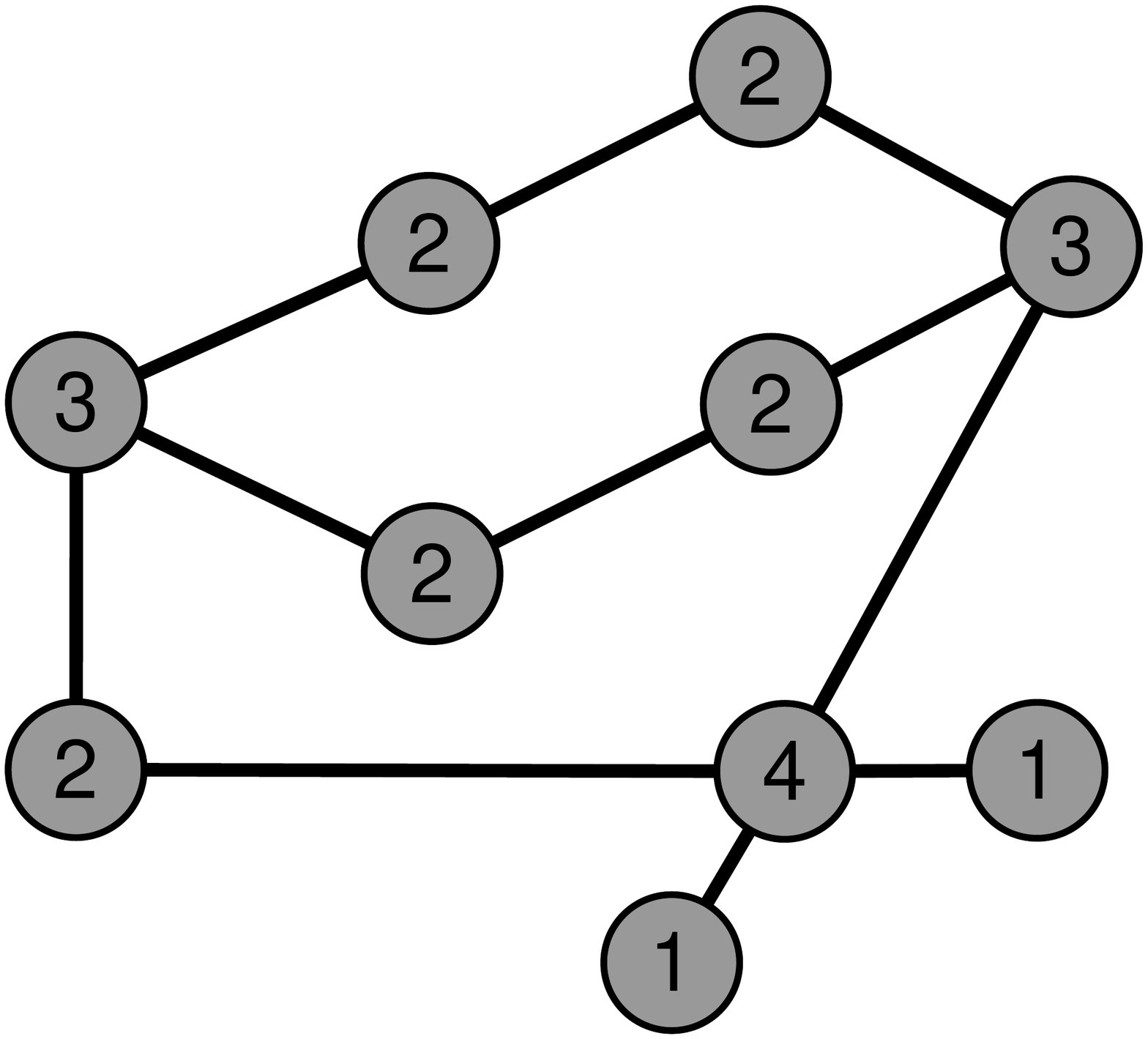}
}
\caption{Representative optimal solutions to Case Studies~\ref{case4} (top) and~\ref{case5} (bottom). The networks at the top are assortative, while the ones at the bottom are disassortative.  The networks to the left feature solutions with isolated nodes.}
\label{fig:case4and5}
\end{figure}

In Case Study~\ref{case4}, we design an assortative network by specifying bounds that are increasing as a function of $q$.
Conversely, in Case Study~\ref{case5}, we construct a disassortative network by specifying bounds that decrease with $q$.

\begin{casestudy}\label{case4}
Let $N=10$ nodes and let $d^L = 0$ and $d^U = N-1=9$ (i.e., no restriction on node degrees).
The set of possible degrees is, thus, $Q = \left\{0,1,\cdots,9\right\}$.
Let also $adn^L_q = 1/2 + \left(2/3\right)q$ and $adn^U_q = 1 + \left(2/3\right)q$, for all $q \in Q^+$, be lower and upper allowable values for the average degree among the neighbors of all nodes that attain degree $q$.
Figure~\ref{fig:case4and5} (top-left) depicts a representative solution to this specification.
We observe that only eight nodes connect in the network, leaving two nodes isolated.
Let us now choose $d^L = 1$, enforcing that no isolated nodes exist.
A representative solution is also shown in Figure~\ref{fig:case4and5} (top-right).
It took less than 0.3 $sec$ CPU to obtain either of these two solutions. \qed
\end{casestudy}

We observe that the solution we obtained for the case of no isolated nodes is in fact a network that is {\em regular of degree 3}; that is, a network where all nodes attain the same degree 3.  This is indeed a valid solution since $\left\{3\right\} \in \left[adn^L_3, adn^U_3\right]=\left[2.5,3.5\right]$. 
The specification for the remaining degree classes, $q \in Q^+ \setminus \left\{3\right\}$, becomes irrelevant in this solution, as these classes were not sampled by any node.  We remark that if more degree diversity is desired, regular graphs can be trivially excluded from consideration by appending the constraint $pd_{1}-pd_{N} \geq n$, where $n \geq 1$ is the minimum span of degrees acceptable for a solution.

\begin{casestudy}\label{case5}
Let $N=10$ nodes and let $d^L = 0$ and $d^U = N-1=9$ (i.e., no restriction on node degrees).
The set of possible degrees is, thus, $Q = \left\{0,1,\cdots,9\right\}$.
Let also $adn^L_q = 4 - \left(2/3\right)q$ and $adn^U_q = 5 - \left(2/3\right)q$, for all $q \in Q^+$, be lower and upper allowable values for the average degree among the neighbors of all nodes that attain degree $q$.
Figure~\ref{fig:case4and5} (bottom-left) depicts a representative solution to this specification.  We observe that only nine nodes connect in the network, leaving one node isolated.
Let us now choose $d^L = 1$, enforcing that no isolated nodes exist.
A representative solution is also shown in Figure~\ref{fig:case4and5} (bottom-right).
It took less than 0.5 $sec$ CPU to obtain either of these two solutions. \qed
\end{casestudy}

We observe that the solutions of Case Study~\ref{case5} do not feature any nodes with degree $7$ or above.  In fact, this was enforced implicitly by the specification, despite the looser setting $d^U = 9$.
Since only connected nodes (i.e., nodes of degree at least $1$) can serve as neighbors of other nodes, the setting $adn^U_{q} < 1$ for all $q \geq 7$ implied that nodes of degree $7$ or above may not have any neighbors, i.e., that no such nodes may exist in the solution.  Additionally, no nodes of degree $6$ exist in these solutions, as the setting $adn^U_{6} = 1$ implies that such nodes may only connect to nodes of degree $1$, forming $6$-star components, which apparently is very restrictive given the rest of the specification (in particular, it is impossible to achieve $adn_{1} \leq adn^U_{1}=4.33$ in a $10$-node network that features a $6$-star component).

We finally remark that, similarly to Case Study~\ref{case2}, in Case Studies~\ref{case4} and~\ref{case5} we were able to obtain solutions that featured disconnected components.  This was possible due to the fact that no explicit restriction was imposed on network properties pertaining to connectivity and/or length of paths between nodes as, for example, was the setting in Case Study~\ref{case3}.

\subsection{Tuning Network Robustness}
Robustness (a.k.a., resilience) of complex networks to attacks or failures is an important subject of study in many fields involving technological networks, most notably in the context of defense applications.
Robustness is usually studied by measuring the effect of removing one or more nodes on the network's functional properties.
In this context, the importance of the individual nodes to the integrity of the network is generally quantified by
using some measure of {\em node centrality}.
The node degrees capture one aspect of node centrality.
However, degree centrality only addresses the importance of the node in its immediate neighborhood.
The {\em closeness centrality} of a node (i.e., the inverse of the average shortest-path length from this node to every other node in the network) can be used as a better metric to measure the node's importance to the network as a whole.
A value of closeness centrality close to 1 indicates that the node is centrally-located, and thus integral to the structural robustness of the network, while a value close to 0 indicates that the node is peripheral to the network, and thus less important from a robustness point of view.
Therefore, by specifying different distributions of closeness centrality for the nodes of a network, we can control the overall network robustness.
In this application, we construct networks with a specified sequence of closeness centralities.

More formally, let $d^L$ and $d^U$, where $0 \leq d^L \leq d^U \leq N-1$, be lower and upper bounds on the node degrees, and let $\mathbf{clc} = \left\{\left[clc_1^L,clc_1^U\right], \left[clc_2^L,clc_2^U\right], \cdots, \left[clc_N^L,clc_N^U\right] \right\}$ be a sequence of ranges for the values of closeness centrality that are to be assumed by the network's nodes.
In other words, we require that there is a one-to-one correspondence between the set of closeness centralities of each node with the $N$ ranges in the sequence $clc$.
The following formulation applies:
\begin{eqnarray}
&\min\limits_{\mathbf{x}} & \sum\limits_{i \in V} \left(sd_i^{-} + sd_i^{+} + siclc_i^{-} + siclc_i^{+}\right) \label{formrob-obj} \\
& \text{s.t.} & \text{Constraints~(\ref{sp-p1})--(\ref{sp-p4}),~(\ref{sp-d1})--(\ref{sp-d4}),~(\ref{sp-o1}),~(\ref{sp-o2})} \;\;\; \forall (i,j) \in E  \nonumber \\
&&\hspace{-0.23cm}\left.\begin{array}{l}
\text{Constraints~(\ref{pd})} \\
\text{Constraints~(\ref{piclc})}\\
\text{Constraints~(\ref{spec-deg})}
\end{array} \;\; \right\} \;\;\;\;\;\;\;\;\;\;\, \forall i \in V  \nonumber \\
&& \text{Constraints~(\ref{spec-clc})} \nonumber\\
&& \text{Constraints~(\ref{symbreak_deg_clc})} \nonumber \\
&& x_{ij} \in \left\{0,1\right\} \;\;\;\;\;\;\;\;\;\;\;\;\;\;\;\;\;\;\;\;\;\;\;\;\;\;\;\; \forall (i,j) \in E \label{formrob-xbinary}\\
&& q_{im} \in \left\{0,1\right\} \;\;\;\;\;\;\;\;\;\;\;\;\;\;\;\;\;\;\;\;\;\;\;\;\;\;\;\, \forall (i,m) \in V \times V \label{formrob-qbinary}\\
&& sd_i^{-},sd_i^{+},siclc_i^{-},siclc_i^{+} \geq 0 \;\;\;\; \forall i \in V \label{formrob-slackpos1}.
\end{eqnarray}

Constraints~(\ref{sp-p1}) through~(\ref{sp-o2}) correspond to the encoding in $w$-variables of shortest-path lengths via the primal/dual formulation analysis presented in Section~\ref{shortestpaths}.
Constraints~(\ref{pd}) and~(\ref{piclc}) define the applicable $p$-variables for the node degrees and the (inverse of) closeness centrality, respectively.
Constraints~(\ref{spec-deg}) impose the specification of overall degree bounds, as they apply to each node independently, while constraints~(\ref{spec-clc}) impose the specification for the closeness centrality sequence.
Constraints~(\ref{symbreak_deg_clc}) eliminate isomorphism.  This is primarily based on node degrees, with nodes of the same degree being arranged in order of decreasing closeness centrality (equivalently, in order of increasing inverse of closeness centrality).
Finally, constraints~(\ref{formrob-xbinary}) and (\ref{formrob-qbinary}) declare the binarity of the $x$- and $q$-variables, while constraints~(\ref{formrob-slackpos1}) declare the non-negativity of the $s$-variables.
The objective~(\ref{formrob-obj}) minimizes the sum of every sum of all slack variables participating in the formulation.

\begin{figure}[!ht]
\captionsetup[subfigure]{labelformat=empty}
\centering
\subfloat[]{
\includegraphics[totalheight=0.32\textheight,trim=0 7.0cm 0 7.0cm,clip=true]{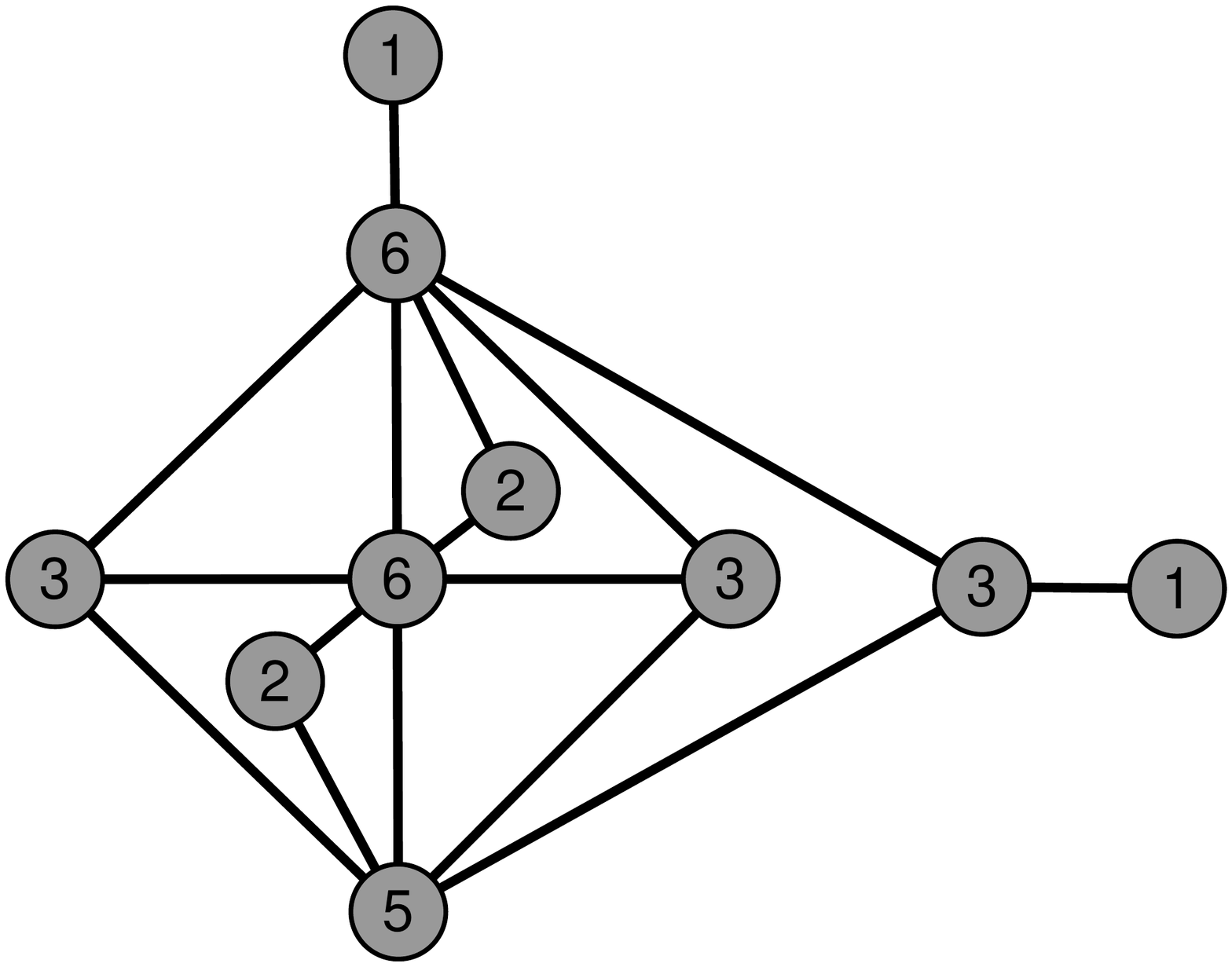}
}
\subfloat[]{
\includegraphics[totalheight=0.27\textheight,trim=0 0.0cm 0 0.0cm,clip=true]{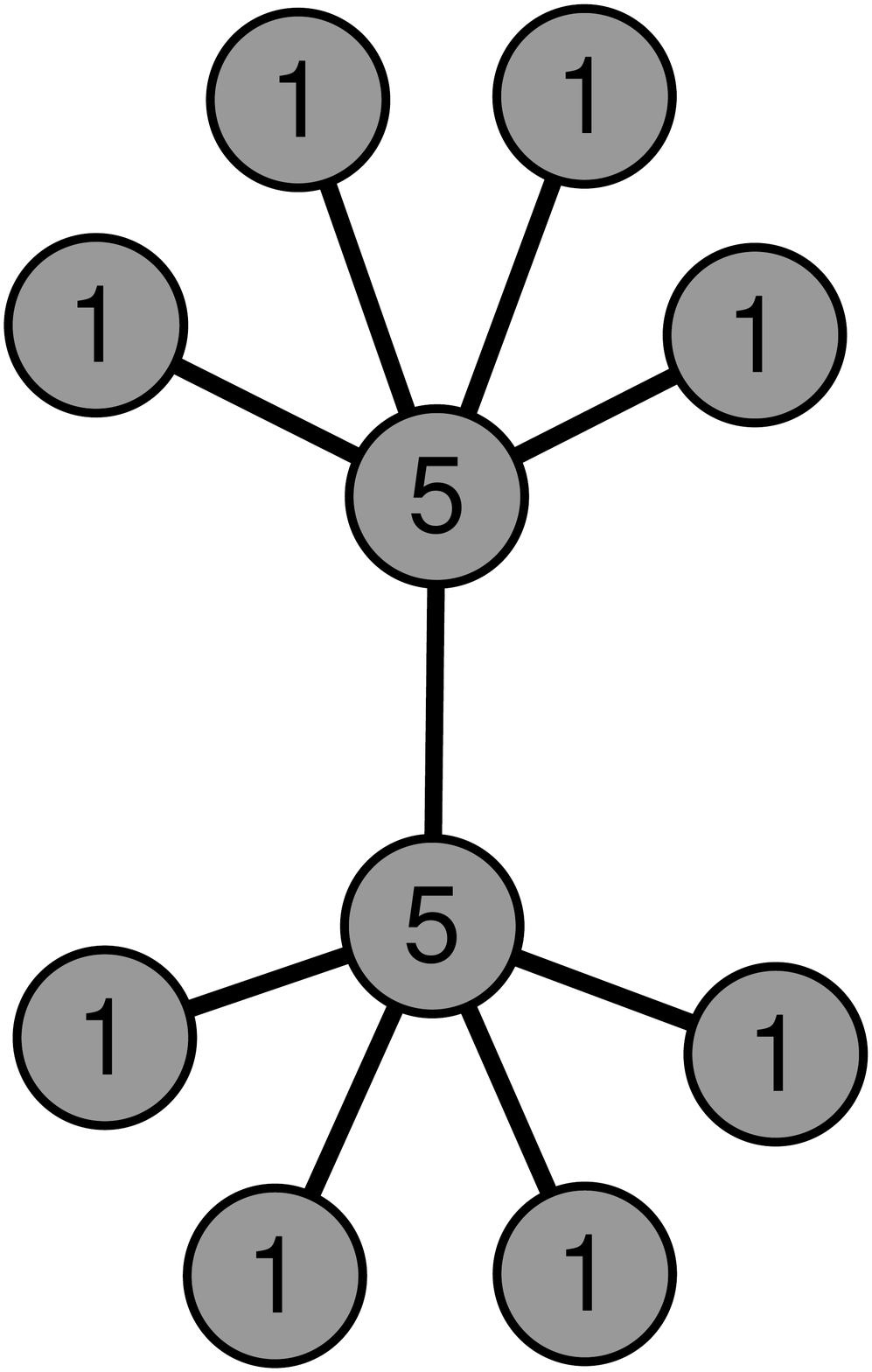}
}
\caption{Representative optimal solutions to Case Studies~\ref{case6} (left) and~\ref{case7} (right).  {\em Left:} The node of degree 6 at the top is the most important for robustness, exhibiting a closeness centrality value of $clc=0.74$, while the node of degree 1 at the right is the least important, with a value of $clc=0.39$. {\em Right:} The two nodes of degree 5 attain a closeness centrality value of $clc = 0.69$ and are distinctively more critical for robustness than the other nodes, all of which attain a much smaller value of $clc = 0.43$.}
\label{fig:case6and7a}
\end{figure}

Case Studies~\ref{case6} and~\ref{case7} demonstrate the use of this formulation.
We aim to design networks that attain two separate sequences for closeness centrality.
The first case study creates a network where the nodes are progressively important for robustness.
The second case study results in a network that possesses a vulnerable link.

\begin{casestudy}\label{case6}
Let $N=10$ nodes and let $d^L = 0$ and $d^U = N-1=9$ (i.e., no restriction on node degrees).
In addition, let the following specification for the closeness centrality sequence, which requires the network nodes to be of increasing importance for robustness:
\begin{equation*}
\left[clc^L_m,clc^U_m \right] = \left[0.35 + m/30, 0.42 + m/30 \right], \;\; \forall m \in V.
\end{equation*}
Figure~\ref{fig:case6and7a} (left) depicts a representative solution for this specification, which was obtained in about 50 $min$ CPU.\footnote{Enabling multi-threaded optimization and allowing the MILP optimizer to utilize all 4 cores of our machine, this number was reduced to less than 4 $min$ CPU.}\qed
\end{casestudy}

\begin{casestudy}\label{case7}
Let $N=10$ nodes and let $d^L = 0$ and $d^U = N-1=9$ (i.e., no restriction on node degrees).
In addition, let the following specification for the closeness centrality sequence, which requires two nodes to be distinctively more important for robustness than the rest of the network:
\begin{equation*}
\left[clc^L_m,clc^U_m \right] = \left\{
\begin{array}{l}
\left[0.10, 0.50 \right], \;\; \forall m \in \left\{1,2,\cdots,8\right\} \\
\left[0.60, 1.00 \right], \;\; \forall m \in \left\{9,10\right\}.
\end{array}\right.
\end{equation*}
Figure~\ref{fig:case6and7a} (right) depicts a representative solution for this specification, which was obtained in 28 $sec$ CPU.\qed
\end{casestudy}

\begin{figure}[!ht]
\captionsetup[subfigure]{labelformat=empty}
\centering
\subfloat[]{
\includegraphics[totalheight=0.27\textheight,trim=0 4.0cm 0 4.0cm,clip=true]{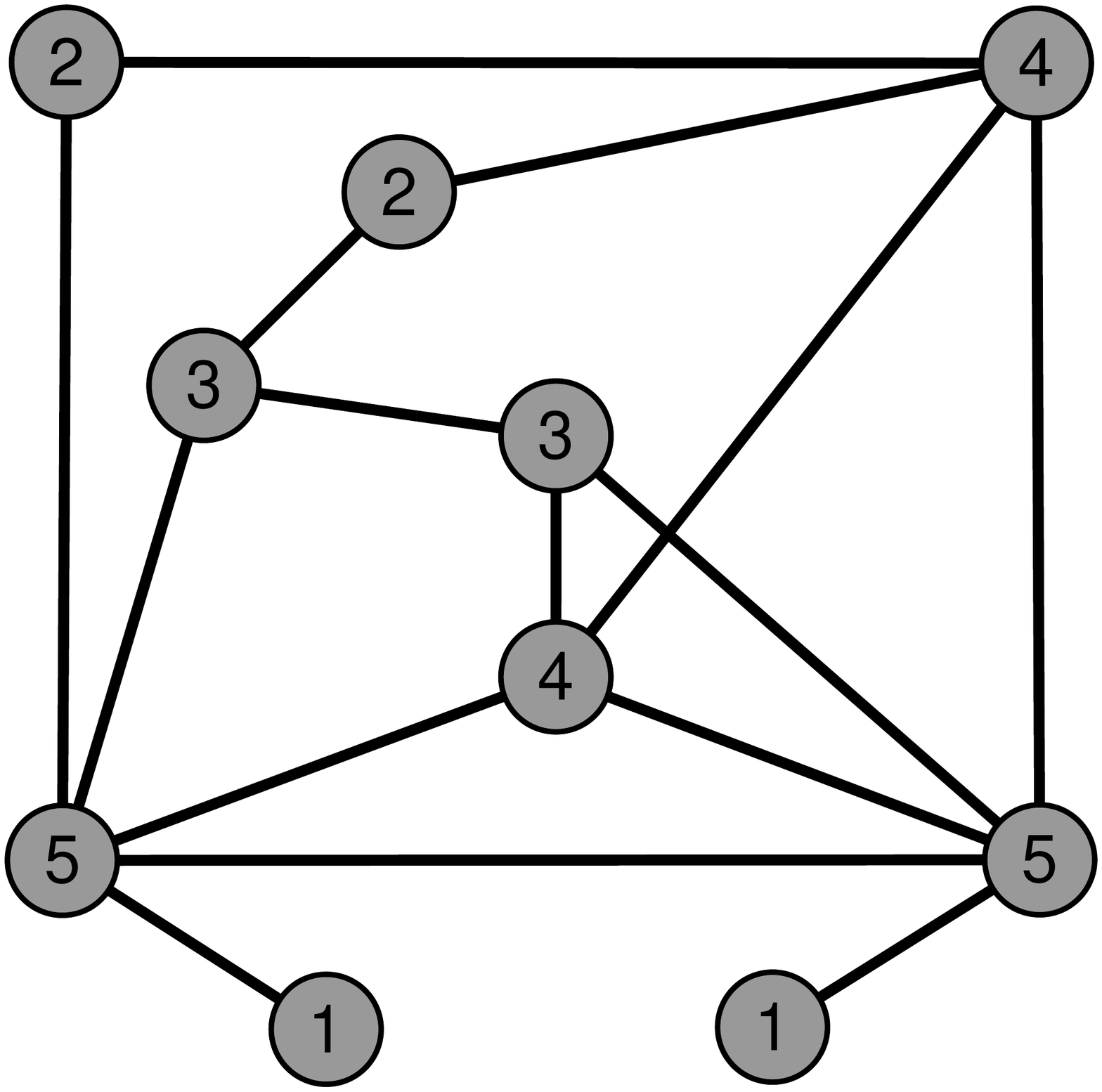}
}
\hspace{2.5cm}
\subfloat[]{
\includegraphics[totalheight=0.27\textheight,trim=0 4.0cm 0 4.0cm,clip=true]{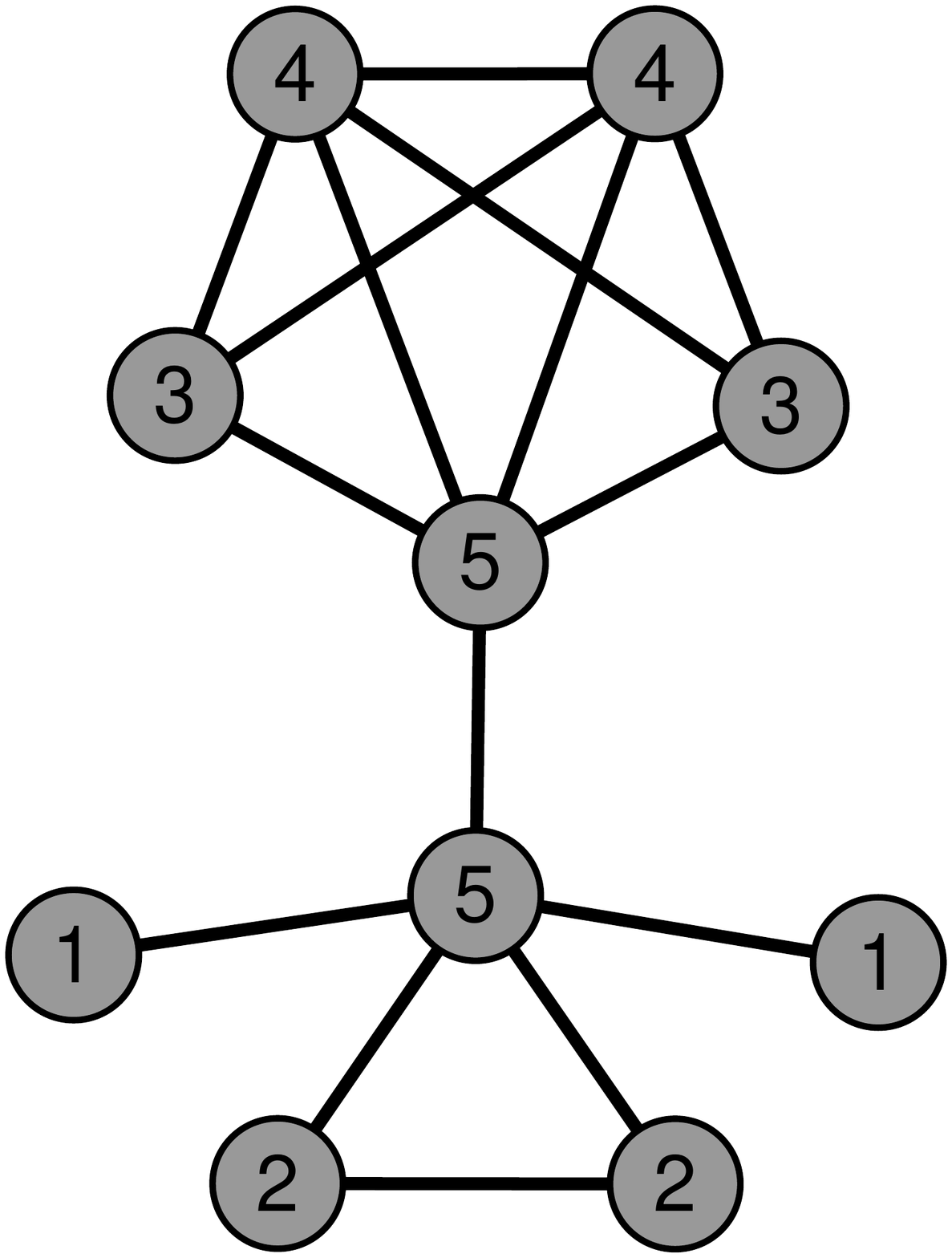}
}
\caption{Optimal solutions to Case Studies~\ref{case6} (left) and~\ref{case7} (right), which additionally possess the degree sequence specified in Case Study~\ref{case3}.}
\label{fig:case6and7b}
\end{figure}

In the above case studies, the degree sequence was not specified.  We will conclude this application by addressing the case where a specific degree sequence is desired in conjunction with the specification of a sequence for the closeness centralities.
Solutions to this problem can be obtained by simply replacing $d^L$ and $d^U$ in constraints~(\ref{spec-deg}) with the desired degree values $d_i$, for all $i \in V$ (equivalently, by replacing  constraints~(\ref{spec-deg}) with constraints~(\ref{spec-degseq})).  We further replace the symmetry-breaking constraints~(\ref{symbreak_deg_clc}) with their simpler version~(\ref{symbreak_deg_clc_fix}).
Figure~\ref{fig:case6and7b} depicts two networks that are also solutions to the specifications of Case Studies~\ref{case6} and~\ref{case7}, but at the same possess the degree sequence specified in Case Study~\ref{case3} (i.e., $\mathbf{d} =\left\{5, 5, 4, 4, 3, 3, 2, 2, 1, 1\right\}$).
These solutions were respectively obtained in about 122 $sec$ and 16 $sec$ CPU.  We note that the former is significantly faster than the corresponding run of Case Study~\ref{case6}, which did not entail a specification for the degree sequence.
Specifying the degrees--and utilizing this fact explicitly in order to break some of the formulation's symmetry--appears in this case to have considerably reduced the feasible search space and to have made the task of identifying a solution easier.

\section{Conclusions}\label{conclusions}
This paper demonstrated the potential and versatility of using Mixed-Integer Linear Optimization to solve the Network Generation Problem, which is the problem of the ``on-demand'' construction of networks that possess specified collective properties.
It was shown how to mathematically express a number of network properties in a format consistent with a MILP optimization formulation,
and how to enforce their specification by restricting them into appropriate values in the context of the optimization formulation.
The particular network properties covered in this paper were chosen due to their popularity among network researchers and were meant to serve as a paradigm.
The modeling techniques that were presented can be utilized in extending the framework to handle a broad class of additional properties that may be relevant to some particular application of interest.

We further compiled complete optimization formulations to specify various levels of connectivity, spread, assortativity and robustness on the networks to be generated, and we solved a number of relevant case studies involving small (10-node) networks to illustrate the overall framework.
Solutions to these case studies were obtained within a few seconds or minutes --depending on the instance-- through the use of uncustomized, ``off-the-shelf,'' optimization software.
There are multiple factors that affect the computational performance of an MILP solution process when dealing with instances of the Network Generation Problem.
These include the size of the network to be constructed, the number and type of variables required to encode the problem (which depends on the application), as well as the restrictiveness of the desired specification.
Future research will be directed towards establishing the limitations of available ``off-the-shelf'' optimization software  and, more importantly, towards developing customized solution methodologies in order to expand the current performance envelope.

\section*{Acknowledgements}\label{funding}
This work was supported by the United States Department of Energy (DE-SC0002097), the Defense Threat Reduction Agency (HDTRA1-07-1-0005) and the National Science Foundation (CBET-0827907).  CEG further acknowledges support from the Berkman Faculty Development Fund.

%
%
\linespread{1}
\small

\bibliographystyle{plain}

\begin{thebibliography}{10}

\bibitem{Aiel00random}
W. Aiello, F. Chung and L. Lu.
\newblock A random graph model for power law graphs.
\newblock {\em Experimental Math}, 10:53--66, 2000.

\bibitem{Albe02statistical}
R. Albert and A.~L. Barab\'{a}si.
\newblock Statistical mechanics of complex networks.
\newblock {\em Reviews of Modern Physics}, 74(1):47--97, 2002.

\bibitem{Aleg86some}
I.~Alegre, M.~A. Fiol and J.~L.~A. Yebra.
\newblock Some large graphs with given degree and diameter.
\newblock {\em Journal of Graph Theory}, 10(2):219--224, 1986.

\bibitem{Badh10impact}
J. Badham and R. Stocker.
\newblock The impact of network clustering and assortativity on epidemic
  behaviour.
\newblock {\em Theoretical Population Biology}, 77(1):71--75, 2010.

\bibitem{Bans09exploring}
S. Bansal, S. Khandelwal and L. Meyers.
\newblock Exploring biological network structure with clustered random
  networks.
\newblock {\em BMC bioinformatics}, 10(1):405, 2009.

\bibitem{Bara99emergence}
A.~L. Barab\'{a}si and R. Albert.
\newblock Emergence of scaling in random networks.
\newblock {\em Science}, 286(5439):509--512, Oct 1999.

\bibitem{Berm84large}
J.C Bermond, C.~Delorme and G.~Farhi.
\newblock Large graphs with given degree and diameter. \{II\}.
\newblock {\em Journal of Combinatorial Theory, Series B}, 36(1):32--48,
  1984.

\bibitem{Blit06sequential}
J. Blitzstein and P. Diaconis.
\newblock A sequential importance sampling algorithm for generating random
  graphs with prescribed degrees.
\newblock {\em Internet Mathematics}, 6(4):489--522, 2011.

\bibitem{Boll80probabilistic}
B.~Bollobas.
\newblock A probabilistic proof of an asymptotic formula for the number of
  labelled regular graphs.
\newblock {\em European J. Combin.}, 1:311–316, 1980.

\bibitem{Bu02distinguishing}
T. Bu and D. Towsley.
\newblock On distinguishing between internet power law topology generators,
  2002.

\bibitem{Calv97modeling}
K. Calvert, M.~B. Doar, A. Nexion, E.~W. Zegura, Georgia Tech,
  and Georgia Tech.
\newblock Modeling internet topology.
\newblock {\em IEEE Communications Magazine}, 35:160--163, 1997.

\bibitem{Chun02connected}
F. Chung and L. Lu.
\newblock Connected components in random graphs with given expected degree
  sequences.
\newblock {\em Annals of Combinatorics}, 6:125--145, 2002.

\bibitem{Deij09random}
M. Deijfen and W. Kets.
\newblock Random intersection graphs with tunable degree distrbution and
  clustering.
\newblock {\em Probability in the Engineering and Informational Sciences},
  23:661--674, 9 2009.

\bibitem{Doro02how}
S.~N. Dorogovtsev, J.~F.~F. Mendes and A.~N. Samukhin.
\newblock {How to construct a correlated net}.
\newblock {\em ArXiv Condensed Matter e-prints}, June 2002.

\bibitem{Falo99power-law}
M. Faloutsos, P. Faloutsos and C. Faloutsos.
\newblock On power-law relationships of the internet topology.
\newblock {\em ACM SIGCOMM Comput. Commun. Rev.}, 29:251--262, 1999.

\bibitem{Fell95large}
M.~Fellows, P.~Hell and K.~Seyffarth.
\newblock Large planar graphs with given diameter and maximum degree.
\newblock {\em Discrete Applied Mathematics}, 61(2):133--153, 1995.

\bibitem{Fu02evolving}
P.~Fu and K.~Liao.
\newblock An evolving scale-free network with large clustering coefficient.
\newblock {\em PNAS}, 99:7821–7826, 2002.

\bibitem{OPTLpaper}
C.~E. Gounaris, K. Rajendran, I.~G. Kevrekidis and
  C.~A. Floudas.
\newblock Generation of networks with prescribed degree-dependent clustering.
\newblock {\em Optimization Letters}, 5(3):435--451, 2011.

\bibitem{Guo06growing}
Q. Guo, T. Zhou, J.-G. Liu, W.-J. Bai, B.-H. Wang and M. Zhao.
\newblock Growing scale-free small-world networks with tunable assortative
  coefficient.
\newblock {\em Physica A},
  371(2):814--822, 2006.

\bibitem{GUROBI}
{Gurobi Optimization}.
\newblock {Gurobi Optimizer 5.6}, 2014.

\bibitem{Haki62realizability}
S.~L. Hakimi.
\newblock On realizability of a set of integers as degrees of the vertices of a
  linear graph. i.
\newblock {\em Journal of the Society for Industrial and Applied Mathematics},
  10(3):496--506, 1962.

\bibitem{Have55remark}
V.~Havel.
\newblock A remark on the existence of finite graphs. (czech).
\newblock {\em Casopis Pest. Mat.}, 80:477–480, 1955.

\bibitem{Holm02growing}
P. Holme and B.~J. Kim.
\newblock Growing scale-free networks with tunable clustering.
\newblock {\em Phys. Rev. E}, 65:026107, 2002.

\bibitem{Hong04factors}
H.~Hong, B.~J. Kim, M.~Y. Choi and H. Park.
\newblock Factors that predict better synchronizability on complex networks.
\newblock {\em Phys. Rev. E}, 69(6):067105, 2004.

\bibitem{Kim04performance}
B.~J. Kim.
\newblock Performance of networks of artificial neurons: the role of
  clustering.
\newblock {\em Phys. Rev. E}, 69(4):045101, 2004.

\bibitem{BTER_2014}
T.~.G. Kolda, A. Pinar, T. Plantenga and C. Seshadhri.  A scalable
  generative graph model with community structure. {\em http://arxiv.org/abs/1302.6636}, 2014.
  
\bibitem{Miao08effects}
Q. Miao, Z. Rong, Y. Tang and J. Fang.
\newblock Effects of degree correlation on the controllability of networks.
\newblock {\em Physica A},
  387(24):6225--6230, 2008.

\bibitem{Mill11efficient}
J.~C. Miller and A. Hagberg.
\newblock Efficient generation of networks with given expected degrees.
\newblock In {\em Algorithms and Models for the Web Graph}, pages 115--126.
  Springer, 2011.

\bibitem{Mil03uniform}
R.~Milo, N.~Kashtan, S.~Itzkovitz, M.~E.~J. Newman and U.~Alon.
\newblock {On the uniform generation of random graphs with prescribed degree
  sequences}.
\newblock {\em eprint arXiv:cond-mat/0312028}, 2003.

\bibitem{nakano}
S.~I. Nakano, R.~Uehara and T. Uno.
\newblock {Efficient Algorithms for a Simple Network Design Problem}.
\newblock {\em Networks}, 62(2):95--104, 2013.

\bibitem{Newm02assortative}
M.~E.~J. Newman.
\newblock Assortative mixing in networks.
\newblock {\em Phys. Rev. Lett.}, 89:208701, 2002.

\bibitem{Newm03mixing}
M.~E.~J. Newman.
\newblock Mixing patterns in networks.
\newblock {\em Phys. Rev. E}, 67:026126, 2003.

\bibitem{Newm03structure}
M.~E.~J. Newman.
\newblock The structure and function of complex networks.
\newblock {\em SIAM Review}, 45(2):167--256, 2003.

\bibitem{Nico12controlling}
V. Nicosia, R. Criado, M. Romance, G. Russo and V. Latora.
\newblock Controlling centrality in complex networks.
\newblock {\em Scientific reports}, 2, 2012.

\bibitem{Past01dynamical}
R.~Pastor-Satorras, A.~Vázquez and A.~Vespignani.
\newblock Dynamical and correlation properties of the internet.
\newblock {\em Phys. Rev. Lett.}, 87(25):258701, 2001.

\bibitem{Pech12influence}
D.~A. Pechenick, J.~L. Payne and J.~H. Moore.
\newblock The influence of assortativity on the robustness of
  signal-integration logic in gene regulatory networks.
\newblock {\em Journal of Theoretical Biology}, 296(0):21--32, 2012.

\bibitem{Pusc08generating}
A. Pusch, S. Weber and M. Porto.
\newblock Generating random networks with given degree-degree correlations and
  degree-dependent clustering.
\newblock {\em Phys. Rev. E}, 77:017101, 2008.

\bibitem{Rava03hierarchical}
E. Ravasz and A.-L. Barabási.
\newblock Hierarchical organization in complex networks.
\newblock {\em Phys. Rev. E}, 67(2):026112, 2003.

\bibitem{Serr05tuning}
M.~A. Serrano and M. Bogun\'a.
\newblock Tuning clustering in random networks with arbitrary degree
  distributions.
\newblock {\em Phys. Rev. E}, 72(3):036133, 2005.

  generative graph model with community structure. {\em http://arxiv.org/abs/1302.6636}, 2014.
  
\bibitem{BTER_2012}
C. Seshadhri, T.~.G. Kolda and A. Pinar. Community structure and scale-free collections of Erd\"{o}s-R\'{e}nyi graphs. {\em Phys. Rev. E}, 85:056109, 2012.

\bibitem{Shir05impacts}
M.~D.~F. Shirley and S.~P. Rushton.
\newblock The impacts of network topology on disease spread.
\newblock {\em Ecological Complexity}, 2(3):287--299, 2005.

\bibitem{Tang02network}
H. Tangmunarunkit, R. Govindan, S. Jamin, S. Shenker and
  Walter Willinger.
\newblock Network topology generators: degree-based vs. structural.
\newblock {\em SIGCOMM Comput. Commun. Rev.}, 32(4):147--159, 2002.

\bibitem{Vige05efficient}
F. Viger and M. Latapy.
\newblock Efficient and simple generation of random simple connected graphs
  with prescribed degree sequence.
\newblock In Lusheng Wang, editor, {\em Computing and Combinatorics}, volume
  3595 of {\em Lecture Notes in Computer Science}, pages 440--449. Springer
  Berlin Heidelberg, 2005.

\bibitem{Volz04random}
E. Volz.
\newblock Random networks with tunable degree distribution and clustering.
\newblock {\em Phys. Rev. E}, 70(5):056115, 2004.

\bibitem{Wang08evolving}
J. Wang, L. Rong and L. Zhang.
\newblock Evolving small-world networks of the local world with tunable
  clustering.
\newblock In {\em Computing, Communication, Control and Management, 2008. CCCM
  '08. ISECS International Colloquium on}, volume~2, pages 369--373, 2008.

\bibitem{Watt98collective}
D.~J. Watts and S.~H. Strogatz.
\newblock Collective dynamics of `small-world' networks.
\newblock {\em Nature}, 393(6684):440--442, 1998.

\bibitem{Webe07generation}
S. Weber and M. Porto.
\newblock Generation of arbitrarily two-point-correlated random networks.
\newblock {\em Phys. Rev. E}, 76:046111, 2007.

\bibitem{Worm80some}
N.~C. Wormald.
\newblock Some problems in the enumeration of labelled graphs.
\newblock {\em Bulletin of the Australian Mathematical Society},
  21(01):159--160, 1980.

\bibitem{Xulv04reshuffling}
R.~Xulvi-Brunet and I.~M. Sokolov.
\newblock Reshuffling scale-free networks: From random to assortative.
\newblock {\em Phys. Rev. E}, 70:066102, 2004.

\bibitem{Zhou10generation}
J. Zhou, G. Xiao, L. Wong, X. Fu, S. Ma and T.~H. Cheng.
\newblock Generation of arbitrary two-point correlated directed networks with
  given modularity.
\newblock {\em Phys. Lett. A}, 374(31):3129--3135, 2010.

\end{thebibliography}

\end{document}